\documentclass{amsart}
\usepackage[foot]{amsaddr}

\title{Optimal reconstruction of general sparse stochastic block models}

\author{Byron Chin$^*$}
\address{$^*$Department of Mathematics, Massachusetts Institute of Technology}
\email{byronc@mit.edu}

\author{Allan Sly$^\dagger$}
\address{$^\dagger$Department of Mathematics, Princeton University}
\email{asly@math.princeton.edu}

\usepackage{arxiv}

\begin{document}

\begin{abstract}%
This paper is motivated by the reconstruction problem on the sparse stochastic block model. Mossel, et. al.  proved that a reconstruction algorithm that recovers an optimal fraction of the communities in the symmetric, 2-community case. The main contribution of their proof is to show that when the signal to noise ratio is sufficiently large, in particular $\lambda^2d > C$, the reconstruction accuracy for a broadcast process on a tree with or without noise on the leaves is asymptotically the same. This paper will generalize their results, including the main step, to a general class of the sparse stochastic block model with any number of communities that are not necessarily symmetric, proving that an algorithm closely related to Belief Propagation recovers an optimal fraction of community labels.
\end{abstract}

\maketitle

\section{Introduction}\label{section: Introduction}

The stochastic block model is one of the most extensively studied random graph models. It extends the well-known Erd\H{o}s--R\'enyi random graph \cite{ER:60} by allowing inhomogeneous structure within the graph. As a result, the block model has facilitated the statistical study of community detection. Given a graph drawn from this distribution with predetermined structure, the goal is to infer the underlying clusters efficiently, see e.g. \cite{HLL:83,SN:97,RCY:11}. Early work in this line focused on the case of fairly dense graphs,  with Dyer and Frieze \cite{DF:89}, Snijders and Nowicki \cite{SN:97}, and Condon and Karp \cite{CK:01} all able to exactly recover the communities of the graph when the average degrees are polynomial in $n$, the number of vertices. A breakthrough of McSherry \cite{M:01} achieved the same performance with only logarithmic degrees. This represents a natural barrier, as the presence of isolated vertices show that exact recovery is not possible with smaller than logarithmic degrees. In fact, standard probabilistic arguments show that when the average degree is constant, a constant fraction of the vertices will be isolated (see e.g. \cite{AS:16}), and as such one cannot hope to recover more than a constant fraction of the communities correctly. 

In this constant degree regime, a natural question to ask is: when is it possible to recover the communities better than a random guess? A detailed answer to this question was conjectured by Decelle, Krzakala, Moore, and Zdeborov\'a \cite{DKMZ:11} based on intuitions from physics.  Part of this conjecture is that above a threshold known as the Kesten--Stigum (KS) bound \cite{KS:66}, it is always possible to detect communities better than random. Moreover, they conjectured that the optimal algorithm is belief propagation. While the analysis of belief propagation remained elusive, it was suggested in \cite{KMMNSLZ:13} that an algorithm based on non-backtracking walks could confirm the conjecture. Following this idea, in a series of independent works by Mossel, Neeman and Sly as well as Massouli\'e, along with Bourdenave and Lelarge, proved this in the case of 2 communities \cite{MNS:15,MNS:18,BLM:15,M:14}. In fact, they also proved that better than random detection is impossible below the KS bound, establishing a threshold in this case. Subsequent works of Abbe and Sandon gave algorithms to perform detection for three or more communities \cite{AS:15a,AS:18}. The recent work \cite{MSS:22} addresses the converse of this question, establishing complex behavior and a number of exact thresholds for 3 and 4 communities.
	
A natural follow up question to this line of work is: when it is possible to detect better than random, what is the best fraction of communities that one can hope to recover? This question of the optimal fraction of communities that can be recovered is studied in \cite{MNS:16} in the simplest case of two symmetric communities. They show that a variant of the belief propagation algorithm recovers an optimal fraction of the vertices in the regime well above the KS bound. Their result is not only one of the first algorithms to guarantee a better than random performance, but it is also the first algorithm that gives a provably optimal recovery. A previous work of the authors \cite{CS:20} generalizes this result to a symmetric block model with any number of communities. 

Since our initial work, a number of new papers have also contributed to this line of research. Notably, Polyanskiy and Yu \cite{YP:22} proved uniqueness of the belief propagation fixed point up to the KS bound for two symmetric communities, fully tightening the result of \cite{MNS:16}. Gu and Polyanskiy \cite{GP:23} proved uniqueness of the fixed point for $q \geq 3$ symmetric communities, extending the results of \cite{CS:20} asymptotically up to the KS threshold. Their techniques differ greatly from ours, relying on a theory of binary and $q$-ary symmetric channels to show that the belief propagation operator is a contraction under a suitably chosen potential function.

\subsection{Our Contribution}
Stemming from the real world motivation for studying sparse versions of the block model, it is common for different communities of nodes to behave differently from one another. For this reason, it makes sense to study block models in which the communities do not behave identically as well. 

In this paper, we further extend the result of \cite{MNS:16}, showing that the same algorithm continues to perform optimally for the block model with any number of potentially asymmetric communities. Our main theorem, which follows from the more general Theorem \ref{thm: main theorem}, is as follows.


\begin{theorem}
Given a community distribution $\pi$ and a stochastic block model with block matrix $dPD_\pi^{-1}$ where $P$ satisfies Conditions \ref{condition: rows far enough apart}, \ref{condition: uniform average degrees}, \ref{condition: no entries of P too small}, \ref{condition: noise is invertible}, \ref{condition: technical} in Section \ref{subsection: conditions}, if all the entries of some non-principal eigenvector  with non-zero eigenvalue are distinct then there exists an algorithm that recovers an asymptotically optimal fraction of community labels whenever $d$ is large enough.
\end{theorem}

To the best of our knowledge, it is the first result proving that an algorithm achieves optimal accuracy in any asymmetric setting. Our result can also be interpreted as supporting evidence for the performance of belief propagation in this greater generality as well. We note here that our asymmetric results do not fully imply our previous results concerning the symmetric block model in \cite{CS:20}, but the techniques are similar and can be synthesized.

Our proof is inspired by the general strategy of \cite{MNS:16}, but the details to carry out each step are much more delicate. In particular, the asymmetry requires an entire matrix of parameters to determine the model, whereas the symmetric block model can be specified by three parameters. This precludes the expression of many quantities of interest using simple, explicit formulas and instead requires the analysis of functions of many variables. Here, we provide a high level overview of their strategy that we follow. 

The proof consists of three main steps. First, we consider the associated broadcast process on a tree. The study of broadcast processes was also initiated by Kesten and Stigum \cite{KS:66}, who derived the threshold for critical behavior on the tree. Follow up works elucidated the connection between the KS bound and the reconstruction problem on the tree, see e.g. \cite{MP:03,EKPS:00,M:04,M:22} for more on the broadcast model. A related problem known as robust reconstruction, considers the question of whether or not information about the root of the tree can be deduced with noise on the leaves. Janson and Mossel showed that this problem is closely related to the same KS bound \cite{JM:04}. In \cite{MNS:16}, it is shown that the probability of recovering the root label correctly in fact tends to the same value for standard reconstruction (without noise) and robust reconstruction. 
This is broken down into two parts, a first and second moment computation followed by a careful analysis of the tree recursion derived from Bayes' rule. 
The final step uses the fact that the neighborhood of a vertex in the sparse stochastic block model can be coupled with a broadcast process on a Galton-Watson tree with high probability. This shows that a rough initial estimate of the communities can be amplified using local applications of belief propagation to give a community prediction with accuracy achieving the upper bound.  

We note that the techniques of \cite{GP:23} for the symmetric block model do not generalize easily to our setting. In their discussion, they mention that there is difficulty in extending their analysis to the case of asymmetric block models, and additional assumptions are required. This supports the presence of the assumptions in our main theorem, to be discussed in Section \ref{subsection: conditions}. 
	
\subsection{Main Definitions}\label{subsection: main definitions}
We make the central definitions that will be used throughout the paper. We start with the basic structure that motivates this paper: the stochastic block model. 

\begin{definition}[Stochastic Block Model]
	The stochastic block model is determined by four parameters: $n \in \mathbb{N}, q \geq 2, \pi $ a probability distribution on $[q]$, and $Q \in [0,1]^{q \times q}$. 
	
	Given a set of vertices $[n] = \{1, \ldots, n\}$, assign them each a community from $[q]$ independently and at random according to a distribution $\pi$. We denote the communities by the vector $\sigma$, where $\sigma_i$ the community of vertex $i$. In particular, for any vertex $i \in [n]$, 
	\[ \pi_j \coloneqq \pi(j) = \Prob{\sigma_i = j}. \]
	Let $D_\pi$ be the $q \times q$ diagonal matrix with $(D_\pi)_{ii} = \pi_i$. 

	Given two vertices in communities $i$ and $j$,  draw the edge between them independently from all other edges with probability according to the edge probability matrix $Q$. In particular,
	\[ Q_{ij} \coloneqq \Prob{uv \in E(G) \middle\vert \sigma_u = i, \sigma_v = j}. \]
	Note that $Q$ must be a symmetric matrix, as we are working with an undirected graph. 
\end{definition}
	
We next need to define a Galton-Watson tree along with the corresponding broadcast process. We follow the formal definition of a tree from the original paper \cite{MNS:16}. The tree is defined as a subset of $\mathbb{N}^*$, the set of finite strings of natural numbers, such that if $s \in T$, then any prefix of $s$ must also be in $T$. Then, for example, the condition that $T$ is $d$-regular can then be expressed as $\forall s \in T$, $|\{ n \in \mathbb{N}: sn \in T \}| = d$. More generally, we can express the number of children of a node $s \in T$ as $\abs{\{ n \in \mathbb{N}: sn \in T \}}$, which will be distributed according to a Poisson random variable. With this formulation, the root is defined as the empty string. 

\begin{definition}[Broadcast Process]
	Let $\rho$ denote this root of the tree, and $\sigma$ be the broadcase labels on the vertices of the tree distributed as follows. Start with the label $\sigma_\rho$ drawn according to $\pi$, and propagate down the tree according to  transition probabilities based on the transition matrix $P$.  
	\[ P_{ij} \coloneqq \Prob{\sigma_v = j \middle\vert \sigma_u = i, uv \in E(T)}  \]
	where $v$ is a child of $u$ in $T$. 
\end{definition}
	
Now that we have defined our structures, we proceed to the quantities related to these structures that we are interested in. First, we define the key quantity of our reconstruction theorem: the Kesten-Stigum threshold. For the Kesten-Stigum threshold, we are interested in the eigenvalues of the matrix $P$. Throughout this paper, we will refer to these eigenvalues as $\lambda_i = \lambda_i(P)$ satisfying 
\[ 1 = \lambda_1 \geq \lambda_2 \geq \cdots \geq \lambda_q \geq 0. \]
The corresponding eigenvectors will be denoted by $\xi_i$ for the eigenvalue $\lambda_i(P)$. The quantity that defines the Kesten-Stigum threshold is $\lambda_2^2d$. 
	
Now that we have a definition of the main quantity $\lambda_2^2d$, we can move to the definitions relevant to the probability of recovering the root. Let $L_k(u)$ be the $k$-th level in the subtree rooted at $u$, and $\sigma_{L_k(u)}$ be the labels 
on this level. 

\begin{definition}[Posterior Probability]
	The posterior probability of the root vertex $u$ being labeled $i$ given the labels on the $m$th level of its subtree is denoted by 
	\[ X_u^{(m)}(i) \coloneqq \Prob{\sigma_u = i \middle \vert \sigma_{L_m(u)}}. \]
\end{definition}

\begin{definition}[MLE]
	The maximum likelihood estimate of the root label given the $m$th sublevel is given by 
	\[ \hat\sigma(m) \coloneqq \arg\max_i X_\rho^{(m)}(i). \]
\end{definition}

\begin{definition}[Reconstruction Probability]
	The probability of obtaining the correct label based on the leaves at distance $m$ is written as  
	\[ E_m \coloneqq \Prob{\hat\sigma(m) = \sigma_\rho} = \E{\max_i X_\rho^{(m)}(i)}. \]Note that this definition applies to any vertex $v \in T$ due to the recursive nature of the tree.
\end{definition} 
	
Now, we move to the setting of robust reconstruction. Denote the noisy labels of the tree by $\tau$. The noisy labels are obtained from the true labels independently and at random. 

\begin{definition}[Noise Matrix]
	The noise matrix is $\Delta = (\Delta_{ij}) \in [0,1]^{q\times q}$ where 
	\[ \Delta_{ij} = \Prob{\tau_u = j \middle\vert \sigma_u = i}. \]
\end{definition}

Relative to these noisy labels, we can similarly define the posterior probability vector $\tilde X_u$, maximum likelihood estimate $\hat\tau$ and probability of correct recovery $\tilde E_m$.

\begin{definition}[Robust Quantities]
	\begin{align*}
		&\tilde X_u^{(m)}(i) \coloneqq \Prob{\sigma_u = i \middle\vert \tau_{L_m(u)}}. \\
		&\hat\tau(m) \coloneqq \arg\max \tilde X_\rho^{(m)}(i). \\
		&\tilde E_m \coloneqq \Prob{\hat\tau(m) = \sigma_\rho} = \E{\max_i \tilde X_\rho^{(m)}}.
	\end{align*}
\end{definition}

At many points in this paper, we will refer to constant coefficients as $C$. In this context, the constant may depend on our fixed parameters $\pi, \delta, \xi, q$. In addition, each time such a constant appears may not refer to the same constant value, even within the same expression. It only refers to some finite expression depending only on the fixed parameters of the model. With these definitions in place, we can now consider necessary conditions for our optimal reconstruction, and then precisely state the main results of the paper. 
	
\subsection{Conditions}\label{subsection: conditions}
We discuss the main conditions we impose on the parameters of the stochastic block model and the motivations and intuition for doing so. 
	
First, we note that it is indeed necessary to impose conditions on the model in order for partial reconstruction to be theoretically possible. Consider, for example, an edge probability matrix $Q$ in which two rows $Q_i$ and $Q_j$ are equal. Then, from the perspective of the graph, the behavior of communities $i$ and $j$ are completely identical. In this case, it will not be possible to distinguish vertices in community $i$ from those in community $j$, and there is a natural obstruction to the accurate reconstruction of labels. This example suggests a natural condition that is necessary to impose: that any two communities should be ``different enough" in distribution. Translated into the language of our parameters, we should have that any two rows of $Q$ be sufficiently far apart. Equivalently, any two rows of $P$ should sufficiently far apart as well. This motivates the following condition, where sufficiently far apart is quantified in terms of $L^1$ distance. 
\begin{condition}\label{condition: rows far enough apart}
	There exists some $\delta > 0$ such that 
	\[ \forall i \neq j, \norm{P_i - P_j}_1 \geq \delta q \lambda_2.  \]
\end{condition}
We foreshadow here that this quantity $\delta$ will have to satisfy a technical inequality depending on some other parameters of the model, such as the initial distribution $\pi$ and the constant $\xi$ defined below. However, this inequality will be loose enough so as not to greatly restrict the class of matrices to which our results apply. Related to this condition, we introduce a notion of distinguishability for the communities that will play a key role in our main theorem.  
\begin{definition}\label{definition: distinguishability}
    We say that the communities can be \textit{$\delta$-distinguished} by a set of vectors $\{v_k\}$ if for every $i$ and $j$, there exists a $k$ such that $\abs{v_k(i) - v_k(j)} > \delta$. 
\end{definition}
It turns out that the above condition implies that the communities can be $\delta$-distinguished by the set of all eigenvectors, see Lemma \ref{lemma: conditions}.
	
Another condition we will impose is that all communities have uniform average degree. At first glance, it may seem that if different communities have differing average degrees, it may be easier to partition vertices into their communities using the information of average degree. However, in the interest of a theoretical optimum, this difference in fact makes the analysis more difficult. This is because of the imbalance of information for some communities relative to others, from the relative lack of edges. More concretely, the assumption of uniform average degrees implies several convenient properties about the corresponding broadcast process on the random tree which will be essential in our analysis. Specifically, we obtain the nice relationship between the matrices $P$ and $Q$ as well as the reversibility of $P$ with respect to the stationary distribution $\pi$.
\begin{condition}\label{condition: uniform average degrees}
	The following quantity is independent of $i$: 
	\[ \sum_{j=1}^q Q_{ij} \pi_j.  \]
\end{condition}
	
The next condition we impose on the matrix $P$ is that there are no entries that are too small. Indeed, consider the worst case scenario that $P_{ij} = 0$ for some $i$ and $j$. Then, consider the setting in which the root vertex of the tree has label $i$, and we are estimating the children of the root. In the small yet positive probability event where one of these children has label $j$, we obtain an overwhelmingly strong signal that the root is in fact not in community $i$. This is a strong obstruction that only occurs with very small probability. In our analysis, this corresponds to a derivative of the recursive function that blows up to infinity, obstructing our step of achieving a uniform bound. We believe the result to be true even without this assumption for the following reason. The tiny probability should balance the terms that explode to infinity in expectation, creating a balance that allows the result to hold. However, a different approach to the analysis will be necessary to circumvent this obstruction. 
\begin{condition}\label{condition: no entries of P too small}
	There exists some $\xi > 0$ such that
	\[ \min P_{ij} > \xi. \]
\end{condition}
	
We ensure that the noise matrix $\Delta$ for robust reconstruction is invertible. This is a natural setting to consider as in most regimes, we require that the noise is in fact small. In this sense, we primarily are concerned with the situation in which $\Delta$ is very close to the identity matrix. Indeed, we will see that this is true as long as the error of our initial estimate is sufficiently small. Thus, we impose the following weaker condition. 
\begin{condition}\label{condition: noise is invertible}
	$\Delta$ is invertible. 
\end{condition}

Finally, the following is a technical condition that will be necessary for our analysis. 
\begin{condition}\label{condition: technical}
    $\frac{2\sqrt{2} \max \pi_k^{3/2} \max \pi_k^{-3/2}}{\xi^3} \lambda_2 < \frac{\delta^2q^2}{8}$.
\end{condition}
	
\subsection{Main Results}\label{subsection: main results}
The goal of this paper is to prove that a repeated local application of Belief Propagation recovers an asymptotically optimal fraction of community labels. Following the ideas of \cite{MNS:16}, who introduced the algorithm, this will be done with three main steps. 
	
\subsubsection*{Step 1: Weighted Majority Method}
The first step is a result in the setting of the random tree. In the symmetric, two-community case, a lower bound on $E_m$ and $\tilde E_m$ was obtained simply by counting the number of leaves in each community and making a guess according to the majority. The main proposition is that with probability close to 1, it is possible to correctly guess the community of the root with this method. 
	
Generalizing to our situation, we require a more sophisticated estimate. Instead of a simple majority of the community counts on the leaves, we compute a weighted sum of these community counts. With the correct weights, which turn out to the be the eigenvectors of our matrix $P$, we can achieve the same result as in the symmetric, two community case. The following proposition concludes the first main step of the proof, and is the culmination of Section \ref{chapter: weighted majority method}. 
\begin{proposition}\label{prop: weighted majority}
	There exists a $C = C(\pi, q, \delta)$ such that if the communities are $\delta$-distinguishable by $\{\xi_i: \lambda_i^2d > C\}$, then there exists $k_0 = k_0(d, \lambda, q)$ so that for all $k \geq k_0$,
	\[ E_k, \tilde E_k \geq 1 - 2qe^{-\frac{\delta^2q^2}{32} \lambda_2^2d}. \]
\end{proposition}
	
\subsubsection*{Step 2: Contraction}
The second main step also takes place on the random tree, and it consists of showing that once the $E_k$ and $\tilde E_k$ are both close to 1, they in fact contract and are equal in the limit as the depth of the tree increases. This concludes the main steps of the proof that take place on the tree, as we will have shown that noisy and non-noisy reconstruction of the root are achievable with the same accuracy. 
	
The idea is to write out the Bayesian recursive formula defining the $X$'s and $\tilde X$'s. Then, we can compute from this formula its derivatives and analyze its size in expectation. By combining the analysis of the derivatives with a Taylor expansion of the recursive function, we will be able to show that on average, we obtain the contraction that we want. This yields the following proposition, which is the culmination of Sections \ref{section: recursion} and \ref{section: contraction}. 
	
\begin{proposition}
	There exists $C_1(\pi, q, \delta, \xi)$ such that if $\frac{2\sqrt{2} \max \pi_k^{3/2} \max \pi_k^{-3/2}}{\xi^3} \lambda_2 < \frac{\delta^2q^2}{8}$, $\lambda_2^2d > C_1$, and the communities are $\delta$-distinguishable by $\{\xi_i: \lambda_i^2d > C\}$ as in Proposition \ref{prop: weighted majority}, then 
	\[ \lim_{n \rightarrow \infty} \E{X_\rho^{(n)}(i) - \tilde X_\rho^{(n)}(i) \middle\vert \sigma_\rho = j} = 0. \]
\end{proposition}
	
\subsubsection*{Step 3: Relation to the Block Model}
Now that the first two steps provide the desired results in the setting of the tree, we only need to relate these results back to the setting of the general block model. This is a standard coupling between the random tree and a tree-like neighborhood of the sparse stochastic block model. Thus, in running our Bayesian belief propagation on a neighborhood of the block model, we will obtain an optimal estimate of the community at which the neighborhood is rooted. Our repeated local applications of belief propagation then recovers an optimal fraction of the communities on our graph drawn from the stochastic block model. 
	
The proof of this follows the strategy from \cite{MNS:15} and \cite{MNS:16}. This is since the coupling of the tree and the graph is a primarily structural result, not depending on the distribution of the communities. Thus, we will only argue the extension from two labels to more labels. This yields the following result: for any vertex in the graph, its size $\Theta(\log n)$ neighborhood along with the true (resp. predicted) labels can be coupled with a Poisson Galton-Watson tree along with a standard (resp. noisy) broadcast process. 
	

	
\subsection{Main Theorem}
Combining the three main steps, we achieve the following theorem, which completes our overarching goal in \cref{section: algorithm}. 
\begin{theorem}\label{thm: main theorem}
	Given $Q, \pi$, if the matrix $P = \frac{1}{d_\pi} QD_\pi$ satisfies Conditions \ref{condition: rows far enough apart}, \ref{condition: uniform average degrees}, \ref{condition: no entries of P too small}, \ref{condition: noise is invertible}, \ref{condition: technical} and the conditions for Sphere Comparison, there then exists $C(\pi, q, \delta)$ and $C_1(\pi, q, \delta, \xi)$ such that if the communities can be $\delta$-distinguished by $\{\xi_i: \lambda_i^2d > C\}$ and $\lambda_2^2d > C_1$, then Algorithm \ref{algorithm} recovers an asymptotically optimal fraction of community labels from a graph $G$ drawn from the stochastic block model with parameters $Q$ and $\pi$. 
\end{theorem}


\begin{remark}
    With the theorem as stated, we have presented it alongside the most fitting black box algorithm: Sphere comparison. However, the Sphere comparison algorithm requires stronger conditions than our contribution. Therefore, our work can be paired with a different black box algorithm with weaker conditions to achieve an analogous result to the above main theorem immediately. 
\end{remark}

\begin{remark}
    As will be discussed later, Condition \ref{condition: technical} can be replaced by an inequality of the form $C_K\lambda_2^K < \frac{\delta^2q^2}{8}$ for any fixed natural number $K$. This is a strictly weaker condition, and the implication of this condition are discussed in the next subsection. 
\end{remark}
	
\subsection{Applications}
We first discuss how our results translate to the setting of symmetric stochastic block models. It is easy to check that the conditions of the main theorem roughly translate to the following. 
\begin{corollary}
	In the symmetric stochastic block model, $G(n, q, \frac{a}{n}, \frac{b}{n})$, if $\lambda_2^2 d > C(q)$, $b \leq a \leq (1 + \frac{1}{10q^2})b$, and $b$ is sufficiently large in terms of $q$, Algorithm \ref{algorithm} recovers an asymptotically optimal fraction of community labels.
\end{corollary}
The key bottleneck for the restraints on the parameters of the block model is the control of the derivatives of the tree recursion discussed in Section \ref{section: recursion}. Combining the analysis from \cite{CS:20}, where the simplicity of the symmetric block model allows for more explicit control of the derivative in the entire parameter range, with our analysis in Section \ref{section: contraction} implies the main theorem for all symmetric block models satisfying $\lambda_2^2d > C(q)$. The new work \cite{GP:23} proves the same result using different techniques. 

Our main contribution is the extension to block models that are not necessarily symmetric, so we exhibit an interesting class of block model matrices that satisfy the conditions of Theorem \ref{thm: main theorem}. In particular, this implies that our result is applicable to a much wider range of matrices than any previous work can guarantee. 

Consider the matrices of the form 
\[ P = \mathbbm{1} \pi^{T} + \frac{1}{\sqrt{d}} \cdot M \]
where $M$ is a matrix whose rows sum to 0, and note that these $P$ arise from edge probability matrices $Q$ of the form 
\[ Q = \frac{1}{n} \left( d \mathbbm{1}\mathbbm{1}^T + \sqrt{d} MD_\pi^{-1} \right).  \]
For $d$ large, we have that the spectrum of $P$ is of the form 1 and eigenvalues of order $O_M(d^{-1/2})$. Moreover, notice that the rows of $P$ are separated by $O_M(d^{-1/2})$ as well. Thus, we can choose $\delta$ in terms of $M$. Then, by choosing $d$ large enough we can satisfy all the conditions specified in the main theorem. If we choose $\pi = \frac{1}{q} \mathbbm{1}$ and $M$ of the form $a\mathbbm{1}\mathbbm{1}^T - aqI$, we recover the case of the symmetric block model with $q$ communities. Relaxing the choices, we get a large family of block models parameterized by $\pi$, $M$, and $d$ to which our theorem applies but previous symmetric results do not apply. Extending this construction further, we can replace the scaling coefficient $\frac{1}{\sqrt{d}}$ by any function $f(d)$ satisfying $f(d) = o(1)$ and $f(d) = \Omega(d^{-1/2})$. This guarantees that as we increase $d$, the quantity $\lambda_2^2d$ will not tend to 0. It only remains to ensure that the matrix $P$ has all entries on the interval $(\xi, 1]$ for some lower bound $\xi$. This is indeed the case when $d$ is large enough, so the scalar $d$ can be chosen large enough so that all the conditions of the theorem are satisfied. 
	
This is just one construction of an edge probability matrix that our new theorem includes. Other examples can be constructed as well, for example with scaling functions $f$ for which $f(d)$ does not necessarily tend to 0. The key takeaway from this demonstration is that the conditions of our theorem are weak enough to accommodate a much larger class of stochastic block model than previously possible. 
	
\section{Preliminaries}\label{chapter: preliminaries}
In this section, we list properties of the broadcast process on the tree that will aid us in future calculations. These facts will be used repeatedly throughout the paper. 
	
	
\begin{fact}[Relationship between $P$ and $Q$]\label{fact: P vs Q}
	Let $d$ be the uniform average degree guaranteed by \cref{condition: uniform average degrees}, and $d_\pi = \sum_{j=1}^q Q_{ij}\pi_j =  \frac{d}{n}$ be the corresponding uniform edge probability. Then \[ P = \frac{1}{d_\pi} QD_\pi. \]
\end{fact}
Since our eventual task is to relate the tree to the block model, this is the one of the most essential relationships. The proof is an application of Bayes' Rule. 
	
	
\begin{fact}[Reversibility]\label{fact:  reversibility}
	$P$ is reversible with respect to $\pi$. 
\end{fact}
This will be helpful in analyzing the broadcast process defined on the tree, as we can then walk both up and down the levels according to the same transition distribution. 
	
\begin{fact}[Diagonalization of $P$]\label{fact: diagonalizing P}
	$P$ is diagonalizable, and moreover there exists a change of basis matrix of the form $UD_\pi^{1/2}$ which suffices, where $U$ is orthogonal. 
\end{fact}
From the Spectral Theorem, we know that $Q$ is diagonalizable by orthogonal matrices. It is natural to desire that $P$ be diagonalizable as well, as we will be working primarily with the eigenvalues of $P$ as opposed to $Q$. This turns out to be the case (see e.g. \cite[Lemma 12.2]{LP:17}).
	
Finally, we show that the rows of $P$ are a perturbation of the stationary distribution $\pi$, on the order of the eigenvalues of $P$. This will be helpful to us in controlling the behavior of certain random variables in terms of $\lambda_2$.
	
\begin{lemma}[Relation of $P$ and $\pi$]\label{lemma: Pclosetopi}
	For all $i, j$,\[ \abs{P_{ij} - \pi_j} \leq \sqrt2 \max \pi_k^{1/2} \max\pi_k^{-1/2} \lambda_2. \]
\end{lemma}
\begin{proof}
	Fact \ref{fact:  reversibility} implies $\pi_j = \sum_{i=1}^q \pi_i P_{ij}$ so we can express our desired quantity as 
	\[ P_{ij} - \pi_j = P_{ij} - \sum_{k=1}^q \pi_k P_{kj} = (1-\pi_i)P_{ij} - \sum_{k\neq i} \pi_kP_{kj} = (vP)_j \]
	where $v = e_i - \pi$ with $e_i$ being the vector of all 0's and a 1 in the $i$th entry. In order to bound the magnitude of this entry, it then suffices to bound the $L^2$ norm of the product vector. This is the quantity 
	\[ \norm{vP}_2 = \norm{vD_\pi^{-1/2}U^T\Lambda UD_\pi^{1/2}}_2. \]
	Here, the diagonal matrix $D_\pi^{1/2}$ can scale the magnitude by at most a factor of $\max_k \pi_k^{-1/2}$. The unitary matrix $U^T$ preserves the $L^2$ norm, and the diagonal matrix $\Lambda$ scales the magnitude by at most $\lambda_2$. Notice here that the constant vector $\mathbbm{1}$ is the eigenvector corresponding to the eigenvalue $		\lambda_1 = 1$. Since $v^T\mathbbm{1} = 0$, the first entry of the product $vD_\pi^{-1/2}U^T$ will always be 0, and thus the contribution from $\lambda_1$ will be removed. Continuing on, we see that again the unitary matrix $U$ preserves $L^2$ norm and the final diagonal matrix $D_\pi^{1/2}$ scales by at most $\max_k \pi_k^{1/2}$. Thus  overall, we have found that 
	\[ \norm{vP}_2 \leq \lambda_2 (\max \pi_k^{1/2})(\max \pi_k^{-1/2}) \norm{v}_2. \]
	We can easily bound that $\norm{v}_2 \leq \sqrt2$. This gives the desired inequality. 
\end{proof}
	
\section{Weighted Majority Method}\label{chapter: weighted majority method}
	
In this section, we carry out the calculation of the lower bound on both $E_m$ and $\tilde E_m$ from a weighted majority estimate. As the probabilities $E_m$ and $\tilde E_m$ are the theoretically optimal probabilities of guessing the root vertex in their respective regimes, any explicit method that provides an accurate guess of the root is a lower bound on the probability.  
	
\subsection{Non-noisy Setting}\label{section: non-noisy setting}
As mentioned in the introduction, we choose the weights of the communities according to an eigenvalue $\xi_i$ of the matrix $P$. We first compute the expectation of such a weighted sum on the $k$th sublevel of the tree. 
	
\begin{lemma}\label{lemma: expectation}
	\[ \E{\sum_{L_k(\rho)} \xi_i^T \sigma_j \middle\vert \sigma_\rho} = \lambda_i^k d^k \xi_i^T\sigma_\rho. \]
\end{lemma}
\begin{proof}
	Using linearity of expectation and the Markovian nature of the broadcast process,
	\begin{align*}
		\E{\sum_{L_k(\rho)} \xi_i^T \sigma_j \middle\vert \sigma_\rho}  &= \E{\abs{L_k(\rho)}} \cdot \E{\xi_i^T\sigma_j \middle\vert \sigma_\rho} = d^k \cdot \xi_i^T\E{\sigma_j \middle\vert \sigma_\rho} = d^k \xi_i^T (P^T)^k \sigma_\rho = d^k \lambda_i^k \xi_i^T\sigma_\rho.
        \end{align*}
\end{proof}
Now that we know the expectation, we want to be able to bound the variance of these sums around their expectation. This will allow us to show that with high probability, the sums corresponding to different root labels live in disjoint intervals. To begin, we compute some statistics about the structure of the tree. The following two lemmas result from inductive computations using the law of total expectation and variance to condition on the children of the root, and taking advantage of the recursive nature of the tree. 
	
\begin{lemma}\label{lemma: level k variance}
	$\Var{\abs{L_{k}(\rho)}} = \sum_{i=k}^{2k-1} d^i$.
\end{lemma}
	
\begin{lemma}\label{lemma: paths between leaves}
	Let $P_{2\ell, k}$ be the number of paths of length $2\ell$ between leaves of the $k$th level of the tree. Then	\[ \E{P_{2\ell, k}} = d^{k + \ell}. \]
\end{lemma}

Armed with the above structural lemmas, we can proceed to the computation of the variance of the weighted sum. 
	
\begin{lemma} \label{lemma: variance}
	\[ \Var{\sum_{L_k(\rho)} \xi_i^T\sigma_j - \lambda_i^kd^k\xi_i^T\sigma_\rho} \leq  \max \pi_k \norm{\xi_i}^2 \cdot d^k \cdot 		\frac{(\lambda_i^2d)^k - 1}{\lambda_i^2d - 1} + o(\lambda_i^{2k}d^{2k}). \]
\end{lemma}
\begin{proof}
	To evaluate this variance, we condition on the structure of the tree, denoted by $T$. Expanding to introduce this conditioning, 
	\begin{align*}
		\Var{\sum_{L_k(\rho)} \xi_i^T\sigma_j - \lambda_i^kd^k\xi_i^T\sigma_\rho} &= \E{\Var{\sum_{L_k(\rho)} 					\xi_i^T\sigma_j - \lambda_i^kd^k\xi_i^T\sigma_\rho \middle\vert T}} + \Var{\E{\sum_{L_k(\rho)} \xi_i^T\sigma_j - \lambda_i^kd^k\xi_i^T\sigma_\rho\middle\vert T}}.
	\end{align*}
	We begin by estimating the second term in the expansion. As before, we can extract the sum according to the size of the $k$th sublevel. Similarly to the proof of Lemma \ref{lemma: expectation}, we obtain
    \[ \E{\sum_{L_k(\rho)} \xi_i^T\sigma_j - \lambda_i^kd^k\xi_i^T\sigma_\rho\middle\vert T} = \abs{L_k(\rho)} \lambda_i^k \xi_i^T\pi - \lambda_i^kd^k \xi_i^T\pi. \]
	Here, the second term is constant, so it will not contribute to the variance. 
    Thus, using Lemma \ref{lemma: level k variance}, we obtain
    \[ \Var{\E{\sum_{L_k(\rho)} \xi_i^T\sigma_j - \lambda_i^kd^k\xi_i^T\sigma_\rho\middle\vert T}} = o(\lambda_i^{2k}d^{2k}). \]
	Moving to the second term, we similarly start off by evaluating the inner variance. 
	\begin{align*}
		\Var{\sum_{L_k(\rho)} \xi_i^T\sigma_j - \lambda_i^kd^k\xi_i^T\sigma_\rho \middle\vert T} &= \Var{\sum_{L_k(\rho)} 	\xi_i^T \sigma_j \middle\vert T} + \Var{d^k\lambda_i^k \xi_i^T \sigma_\rho \middle\vert T} - 2\Cov{\sum_{L_k(\rho)} \xi_i^T\sigma_j}{d^k\lambda_i^k \xi_i^T\sigma_\rho \middle\vert T}. \\
		\intertext{Extracting the structure of the tree, we get}
		&= \sum_{L_k(\rho)} \Var{\xi_i^T\sigma_j} + \sum_{j \neq j'} \Cov{\xi_i^T \sigma_j}{\xi_i^T\sigma_{j'}} \\
		&\qquad + d^{2k}\lambda_i^{2k} \Var{\xi_i^T\sigma_\rho} - 2 \sum_{L_k(\rho)} d^k\lambda_i^k \Cov{\xi_i^T\sigma_j}			{\xi_i^T\sigma_\rho}.
	\end{align*}
	For the variance terms, we may simply upper bound by the expectation of the square. In the case of \newline $\Cov{\xi_i^T\sigma_j}{\xi_i^T\sigma_\rho}$, we express the distribution of $\sigma_j$ in terms of $\sigma_\rho$ to obtain a factor of $\lambda_i^k$ and $\Var{\xi_i^T\sigma_\rho}$ which we upper bound as before. Finally, with the term $\Cov{\xi_i^T \sigma_j}{\xi_i^T\sigma_{j'}}$, we take advantage of the reversibility of the broadcast process. We can compute this term equivalently by starting at the root $\sigma_j$ and walking $d(j, j')$ steps to the leaf vertex $j'$. By the previous reasoning, this yields a factor of $\lambda_i^{d(j, j')}$ along with an expectation of $					(\xi_i^T\sigma_j)^2$. Combining these estimates, we get
	\begin{align*}
		\Var{\sum_{L_k(\rho)} \xi_i^T\sigma_j - \lambda_i^kd^k\xi_i^T\sigma_\rho \middle\vert T} &\leq \abs{L_k(\rho)} \E{(\xi_i^T\sigma_j)^2} + \sum_{j \neq j'} \lambda_i^{d(j, j')} \E{(\xi_i^T\sigma_j)^2} \\
		&\qquad+ d^{2k} \lambda_i^{2k} \E{(\xi_i^T\sigma_\rho)^2} - 2\abs{L_k(\rho)}d^{k}\lambda_i^{2k} 						\E{(\xi_i^T\sigma_\rho)^2} \\
		&= \abs{L_k(\rho)} \E{(\xi_i^T\sigma_j)^2} + \sum_{\ell=1}^k \lambda^{2\ell} P_{2\ell, k} \E{(\xi_i^T \sigma_j)^2}  \\
		&\qquad+ d^{2k} \lambda_i^{2k} \E{(\xi_i^T\sigma_\rho)^2} - 2\abs{L_k(\rho)}d^{k} \lambda_i^{2k} \cdot \E{(\xi_i^T 		\sigma_\rho)^2} 
	\end{align*}
	where $P_{2\ell, k}$ is the number of paths of length $2\ell$ between leaves of the $k$th level of the tree. 
    Using Lemma \ref{lemma: paths between leaves}, we obtain that 
	\begin{align*}
		\E{\Var{\sum_{L_k(\rho)} \xi_i^T\sigma_j - \lambda_i^kd^k\xi_i^T\sigma_\rho \middle\vert T}} &= 						d^k\E{(\xi_i^T\sigma_j)^2} + \sum_{\ell=1}^k \lambda^{2\ell} d^{\ell+k} \E{(\xi_i^T \sigma_j)^2}  \\
		&\qquad+ d^{2k} \lambda_i^{2k} \E{(\xi_i^T\sigma_\rho)^2} - 2d^{2k} \lambda_i^{2k} \cdot \E{(\xi_i^T \sigma_\rho)^2} \\
		&= \E{(\xi_i^T \sigma_\rho)^2}\cdot d^k \left[ 1 + \sum_{\ell=1}^k (\lambda_i^2d)^\ell - d^k\lambda_i^{2k}\right] \\
            &= \left(\sum_{k=1}^q \pi_k\xi_i(k)^2\right) \cdot d^k \left[ \frac{(\lambda_i^2d)^{k+1} - 1}{\lambda_i^2d - 1} - (\lambda_i^2d)^k \right] \\
		&= \left(\sum_{k=1}^q \pi_k\xi_i(k)^2\right) \cdot d^k \cdot \frac{(\lambda_i^2d)^k - 1}{\lambda_i^2d - 1} \\
		&\leq \max \pi_k \norm{\xi_i}^2 \cdot d^k \cdot \frac{(\lambda_i^2d)^k - 1}{\lambda_i^2d - 1}.
	\end{align*}
	Combining the two estimates derived above gives the desired bound.
\end{proof}
	
	\subsection{Noisy Setting}
	In the setting where we only know the noisy labels on the leaves, we derive similar results about the expectation and variance of the weighted sums. In this case, we will use different weights in order to achieve the same result. As we will see from the computation, it is natural to consider the weights from $(\Delta^T)^{-1}\xi_i$, the eigenvectors transformed by the inverse of the noise matrix. This is also where we use Condition \ref{condition: noise is invertible} that the noise matrix be invertible. 
	\begin{lemma}\label{lemma: noisy expectation}
		\[ \E{\sum_{L_k(\rho)} ((\Delta^T)^{-1}\xi_i)^T \tau_j \middle\vert \sigma_\rho} = \lambda_i^k d^k \xi_i^T\sigma_\rho. \]
	\end{lemma}
	\begin{proof}
		Calculating as in the non-noisy case with the new weights, we see that 
		\begin{align*}
			\E{\sum_{L_k(\rho)} \xi_i^T \tau_j \middle\vert \sigma_\rho}  &= \E{\abs{L_k(\rho)}} \cdot \E{((\Delta^T)^{-1}\xi_i)^T \tau_j \middle\vert \sigma_\rho} = d^k \cdot ((\Delta^T)^{-1}\xi_i)^T\E{\tau_j \middle\vert \sigma_\rho} = d^k \xi_i^T \Delta^{-1}\Delta(P^T)^k \sigma_\rho.
		\end{align*}
		Here, the noise from the label and the inverse noise we introduced to the weights cancel, and we recover the original expectation.
	\end{proof}
	
	\begin{lemma} \label{lemma: noisy variance}
		\[ \Var{\sum_{L_k(\rho)} ((\Delta^T)^{-1}\xi_i)^T \tau_j - \lambda_i^kd^k\xi_i^T\sigma_\rho} \leq \norm{\Delta^{-1}}^2 \norm{\xi_i}^2 \cdot d^k \cdot \frac{(\lambda_i^2d)^k - 1}{\lambda_i^2d - 1} + o(\lambda_i^{2k}d^{2k}). \]
	\end{lemma}
	\begin{proof}
		The proof is exactly the same as the proof of \cref{lemma: variance}, except where we have $\E{(\xi_i^T\sigma_j)^2}$ we instead have $\E{(((\Delta^T)^{-1}\xi_i)^T\tau_j)^2}$. This evaluates to 
		\begin{align*}
			\E{(((\Delta^T)^{-1}\xi_i)^T\tau_j)^2} &= \sum_{k=1}^q (\Delta\pi)_k ((\Delta^T)^{-1}\xi_i)_k^2 \leq \sum_{k=1}^q ((\Delta^T)^{-1}\xi_i)_k^2 = \norm{(\Delta^T)^{-1}\xi_i}^2.
		\end{align*}
        Substituting this bound into the calculation above gives the inequality.  
	\end{proof}
	
	\subsection{Weighted Majority Estimate}
	Now that we have the exact locations of the mean weighted sum and a bound on the variance around these weights, we can apply a second moment argument to estimate the probability that these sums lie on disjoint intervals. This corresponds to the weighted sum on the leaves lying close to a unique expectation corresponding to a specific value of $\ell$, for which we can then say that with high probability, the root indeed had label $\ell$. 
	
	We start with a lemma relating our imposed condition, \cref{condition: rows far enough apart}, to properties of the eigenvectors of $P$. Since the expectations computed from \cref{lemma: expectation} and \cref{lemma: noisy expectation} are based on entries of these eigenvectors, the following lemma tells us that these expectations are well-spaced for distinct values of $\sigma_\rho$. 
	\begin{lemma}\label{lemma: conditions}
		Our condition that $\forall \ i \neq  j, \ \norm{P_i - P_j} > \delta q\lambda_2$ implies $\norm{\xi_\ell} \leq \max \pi_k^{-1/2}$ and \\  $(\xi_\ell(i) - \xi_\ell(j))^2 > \delta^2$ for some $\ell$. In particular, the communities can be $\delta$-distinguished by the set of all eigenvectors. 
	\end{lemma}
	\begin{proof}
		Writing the entries of $P$ in terms of the diagonalization from \cref{fact: diagonalizing P}, we get
		\begin{align*}
			\norm{P_i - P_j}_1 &= \sum_{k=1}^q \abs{P_{ik} - P_{jk}} = \sum_{k=1}^q \abs{\sum_{\ell=1}^q \lambda_\ell \left( \frac{\pi_k}{\pi_i}\right)^{1/2} u_\ell(i) u_\ell(k) - \sum_{\ell=1}^q \lambda_\ell \left( \frac{\pi_k}{\pi_j}\right)^{1/2} u_\ell(j)u_\ell(k)} \\
			&= \sum_{k=1}^q \abs{\sum_{\ell=1}^q \lambda_\ell u_\ell(k) \pi_k^{1/2} \left[ \pi_i^{-1/2} u_\ell(i) - \pi_j^{-1/2} u_\ell(j) \right]}.
		\end{align*}
		Noting that the entries $\pi_i^{-1/2}u_\ell(i)$ are precisely the entries of the eigenvalues $\xi_\ell$. Thus,
		\begin{align*}
			\norm{P_i - P_j}_1 &= \sum_{k=1}^q \abs{\sum_{\ell=1}^q \lambda_\ell u_\ell(k) \pi_k^{1/2}(\xi_\ell(i) - \xi_\ell(j))} \leq \sum_{k, \ell = 1}^q \lambda_\ell u_\ell(k) \pi_k^{1/2} \abs{\xi_\ell(i) - \xi_\ell(j)}.
		\end{align*}
		Note here that all terms involving the largest eigenvalue $\lambda_1 = 1$ are canceled out, as the associated eigenvector is constant, and so the term $\xi_1(i) - \xi_1(j)$ will always be 0. Note also that 
		\[ \sum_{k=1}^q u_\ell(k)\pi_k^{1/2} = u_\ell^T \pi^{1/2} \leq \norm{u_\ell}_2\norm{\pi^{1/2}}_2 = 1. \]
		Now, if for all $\ell$, we had that $\abs{\xi_\ell(i)  - \xi_\ell(j)} < \delta$, then we would have $\norm{P_i - P_j}_1 < \delta\sum_{\ell=2}^q \lambda_\ell \leq\delta q\lambda_2$. Thus, we can deduce that if $\norm{P_i - P_j} > \delta q\lambda_2$ then there must exist some eigenvector $\xi_\ell$ such that $\abs{\xi_\ell(i) - \xi_\ell(j)} \geq \delta$. Moreover, the eigenvectors obtained from our diagonalization have the prescribed norm. 
	\end{proof}
	
	\begin{proposition}\label{prop: majority}
		There exists an $\ell$ such that
		\[ \liminf_{k \rightarrow \infty}E_k >  1 - \frac{4\max \pi_k^{1/2}}{\delta^2}\frac{1}{\lambda_\ell^2d-1}. \]
		In particular $\ell = \arg\min\{|\lambda_i|\}_{i \in I}$ where the communities can be $\delta$-distinguished by $\{\xi_i\}_{i \in I}$. The same holds for $\tilde E_k$. 
	\end{proposition}
	\begin{proof}
		By the calculation in \cref{lemma: expectation}, given two distinct labels of the root $\ell_1$ and $\ell_2$, for any eigenvector $\xi_i$ the mean differs by $d^k \lambda_i^k (\xi_i(\ell_1) - \xi_i(\ell_2))$. We want the variance around each mean to be small enough so that with high probability the obtained values with different root labels will be separated. In particular, we want the probability that the sum is concentrated within half this difference around the mean. The probability of overlap is then bounded using Chebyshev's inequality by 
		\begin{align*}
			\frac{\max \pi_k \norm{\xi_i}^2 d^k \cdot \frac{(\lambda_i^2d)^k-1}{\lambda_i^2d-1}+o(\lambda_i^{2k}d^{2k}) }{\left( \frac{1}{2} (\lambda_id)^k(\xi_i(\ell_1) - \xi_i(\ell_2))\right)^2} \xrightarrow{k \rightarrow \infty} \frac{4\max \pi_k\norm{\xi_i}^2}{(\xi_i(\ell_1) - \xi_i(\ell_2))^2)} \cdot \frac{1}{\lambda_i^2d-1}.
		\end{align*}
		In order for this to be small as we desired, we need to satisfy a couple conditions. In particular, we require that 
		\[ \forall \ell_1 \neq \ell_2 \exists i \ s.t. \norm{\xi_i} < c, \ \xi_i(\ell_1) - \xi_i(\ell_2) > \delta,\  \lambda_i^2d > C.  \]
		By \cref{lemma: conditions} have that $\norm{\xi_i} \leq \max \pi_k^{-1/2}$, $(\xi_i(\ell_1) - \xi_i(\ell_2))^2 > \delta^2$ for some $i$. Substituting in these bounds gives the result. By the exact same argument, we get the result for the noisy version, $\tilde E_k$ as well. 
	\end{proof}
	In particular, this implies that observing the distribution of the leaves in relation to the eigenvectors $\xi_i$, we will be able to distinguish between the community of the root being $\ell_1$ and $\ell_2$ with probability at least $1 - O((\lambda_i^2d)^{-1})$. Since this holds for any two communities, we can identify the true community by looking at the weighted sums according to each of these eigenvectors. 
	
	\subsection{Iterated Estimate}
	From the first lower bound we obtain from Chebyshev's Inequality in \cref{prop: majority}, we have a probability of error that is polynomially small. However, we would like to amplify this lower bound to achieve an exponentially small probability of error. The idea is to first use the initial bound to control the guesses on the children of the root. Then, using these predictions, bound the probability of error in reconstructing the root using Hoeffding's Inequality. The first lemma we prove towards this goal is a result about the distribution of predictions on the children of the root. In particular, we show that the distribution of our guesses on the children is very close to the true distribution given by the rows of $P$. 
	\begin{lemma}\label{lemma: numChildCancel}
		Suppose that $\sigma_\rho = i$. There exists an $\ell$ such that for layer $L$ deep enough and any child $v$, 
		\[ \abs{\E{\Prob{\sigma_v = j \middle\vert \sigma_L^{(i)}}} - P_{ij}} \leq \frac{8q\sqrt{2}\max\pi_k^{3/2} \max\pi_k^{-1}}{\delta^2} \cdot \lambda_2 \cdot \frac{1}{\lambda_\ell^2d}. \]
	\end{lemma}
	\begin{proof}
		From \cref{prop: majority}, we know that for some index $\ell$ and for layer $L$ deep enough, we have that 
		\[ \Prob{\sigma_\rho = j \middle\vert \sigma_L^{(j)}} = 1 - \epsilon \geq  1 - \frac{8\max \pi_k^{1/2}}{\delta^2} \cdot \frac{1}{\lambda_\ell^2d}. \]
		Note that this implies moreover that 
		\[ \Prob{\sigma_\rho = k \middle\vert \sigma_L^{(j)}} \leq \epsilon \leq \frac{8\max \pi_k^{1/2}}{\delta^2} \cdot \frac{1}{\lambda_\ell^2d}, \quad \forall k \neq j. \]
		We can then compute the desired expectation. 
		\begin{align*}
			\E{\Prob{\sigma_v = j \middle\vert \sigma_{L(\rho)}^{(i)}}} &= \E{\sum_{k=1}^q P_{ik} \Prob{\sigma_v = j \middle\vert \sigma_{L-1(v)}^{(k)}}}. \\
			\intertext{Introducing a term defined by the stationary distribution $\pi$ as opposed to the distribution $P_i$, }
			&= \E{\sum_{k=1}^q \pi_k \Prob{\sigma_v = j \middle\vert \sigma_{L-1}^{(k)}} + \sum_{k=1}^q (P_{ik}-\pi_k) \Prob{\sigma_v = j \middle\vert \sigma_{L-1}^{(k)}}} \\
			&=  \E{\sum_{k=1}^q \pi_k \Prob{\sigma_v = j \middle\vert \sigma_{L-1}^{(k)}}} + \E{ \sum_{k=1}^q (P_{ik}-\pi_k) \Prob{\sigma_v = j \middle\vert \sigma_{L-1}^{(k)}}}. \\
			\intertext{Notice that in the first term, we are averaging a probabililty conditional on $\sigma_\rho = k$ over the initial distribution of $\sigma_\rho$. Thus, in expectation we simply get the expectation without the prior conditioning.}
			&= \E{\Prob{\sigma_v = j \middle\vert \sigma_{L-1}}} + \sum_{k=1}^q (P_{ik} - \pi_k)\E{\Prob{\sigma_v = j \middle\vert \sigma_{L-1}^{(k)}}} \\
			&= \pi_j + (P_{ij} - \pi_j)\E{\Prob{\sigma_v = j \middle\vert \sigma_{L-1}^{(j)}}} + \sum_{k \neq j} (P_{ik} - \pi_k)\E{\Prob{\sigma_v = j \middle\vert \sigma_{L-1}^{(k)}}}. \\
			\intertext{Since we assumed that the expectation in the middle term is equal to $1-\epsilon$, we can substitute this in to get}
			&= P_{ij} - (P_{ij} - \pi_j) \epsilon +  \sum_{k \neq j} (P_{ik} - \pi_k)\E{\Prob{\sigma_v = j \middle\vert \sigma_{L-1}^{(k)}}}.
		\end{align*}
		Now that we have introduced the term $P_{ij}$, we can find the desired difference and calculate that 
		\begin{align*}
			\abs{\E{\Prob{\sigma_v = j \middle\vert \sigma_L^{(i)}}}  - P_{ij}} &= \abs{- (P_{ij} - \pi_j) \epsilon +  \sum_{k \neq j} (P_{ik} - \pi_k)\E{\Prob{\sigma_v = j \middle\vert \sigma_{L-1}^{(k)}}}} \\
			&\leq \abs{P_{ij} - \pi_j} \epsilon + \sum_{k \neq j} \abs{P_{ik} - \pi_k} \E{\Prob{\sigma_v = j \middle\vert \sigma_{L-1}^{(k)}}}.
		\end{align*}
		By \cref{lemma: Pclosetopi} and the bound on $\epsilon$ from \cref{prop: majority},
		\begin{align*}
			\abs{\E{\Prob{\sigma_v = j \middle\vert \sigma_L^{(i)}}}  - P_{ij}} &\leq q \sqrt2 \max \pi_k^{1/2} \max\pi_k^{1/2} \lambda_2 \cdot \frac{8\max \pi_k^{-1/2}}{\delta^2} \cdot \frac{1}{\lambda_\ell^2d} \\
			&= \frac{8q\sqrt{2}\max\pi_k^{3/2} \max\pi_k^{-1}}{\delta^2} \cdot \lambda_2 \cdot \frac{1}{\lambda_\ell^2d}
		\end{align*}
		which is the desired upper bound. 
	\end{proof}
	
	Now, we can use this estimate on the distribution of the estimates on the children to carry out our iterated majority estimate. 
	\begin{proposition} \label{prop: doubleMajority}
		There exists a $C = C(\pi, q, \delta)$ such that if the communities are $\delta$-distinguishable by $\{\xi_i: \lambda_i^2d > C\}$, then there exists $k_0 = k_0(d, \lambda, q)$ so that for all $k \geq k_0$,
		\[ E_k, \tilde E_k \geq 1 - 2qe^{-\frac{\delta^2q^2}{32} \lambda_2^2d}. \]
	\end{proposition}
	\begin{proof}
		We use the above majority estimate to provide a tight bound on the optimal probability of correctly guessing the root label as 1. The idea of this algorithm is to guess the labels on each of the roots children. Then, based on the distribution on these children, guess the label of the root. We have the lower bound on the optimal guessing from above, and we can then compute the probability of each child being guessed as a particular label.
		
		We guess the roots with the following process. First, from observing the leaves of the tree, we calculate the conditional probabilities that the children is in each of the $q$ communities. Then, we draw a guess for this vertex based on the conditional distribution. From \cref{lemma: numChildCancel} above, we know that from this process, given that the root is in community $i$, there is an index $\ell$ such that the probability a given child is guessed to be in community $j$ is 
		\[ \tilde P_{ij} \coloneqq \E{\Prob{\sigma_v = j \middle\vert \sigma_L^{(i)}}} = P_{ij} + O(\lambda_2 \cdot \frac{1}{\lambda_\ell^2d}). \]
		In particular, the number of children we guess to be in community $j$ will be distributed like $\text{Bin}(d,\tilde P_{ij})$. Recall that we assumed the distributions $P_i$ differ in $L^1$ norm by at least $\delta q \lambda_2$. Then, as long as $\lambda_\ell^2d$ is large enough, we have that the distributions of guesses also differ by $L^1$ norm at least $\frac{\delta q}{2} \cdot \lambda_2$. In particular, this requires us to have
		\[ \frac{8q\sqrt{2}\max\pi_k^{3/2} \max\pi_k^{-1}}{\delta^2} \cdot \lambda_2 \cdot \frac{1}{\lambda_\ell^2d} < \frac{\delta q}{4} \lambda_2 \iff \lambda_\ell^2d > \frac{32\sqrt{2}\max \pi_k^{3/2}\max\pi_k^{-1}}{\delta^3}. \]
        Notice that we need this to hold for a $\delta$-distinguishable set of eigenvectors $\xi_\ell$, which is exactly the assumption with $C = \frac{32\sqrt{2}\max \pi_k^{3/2}\max\pi_k^{-1}}{\delta^3}$. 
		
		Let $D$ be the number of children of the root and suppose first that this is fixed. Now, let $N_i \sim \text{Bin}(D, \tilde P_{ij})$. Applying Hoeffding's inequality, we know that
		\[ \Prob{\abs{N_i - D\tilde P_{ij}} > \alpha\lambda_2 D} = \Prob{\abs{N_i - \mathbb{E}[N_i]} > \alpha\lambda_2 D} < 2e^{-2(\alpha\lambda_2 D)^2 / D} = 2e^{-2\alpha^2 \lambda_2^2D}. \] 
		Applying the union bound, we get that 
		\[ \Prob{\bigcup_{i=1}^q \left\lbrace\abs{N_i - D\tilde P_{ij}} > \alpha\lambda_2 D\right\rbrace } \leq \sum_{i=1}^q \Prob{\abs{N_i - Dp_i} > \alpha\lambda_2 D} < \sum_{i=1}^q 2e^{-2\alpha^2\lambda_2^2 D} = 2qe^{-2\alpha^2\lambda_2^2D}. \]
		Then, we just need to choose $\alpha$ small enough so that a discrepancy of $\alpha \lambda_2 D$ on the distributions conditioned on differing communities of the root will not overlap. From our choices above, it suffices to choose $\alpha = \frac{\delta q}{4}$. Thus, the probability we can correctly recover the root with this method conditioned on $D$ is at least $1-2qe^{-\frac{\delta^2q^2}{8} \lambda_2^2D}$. 
		
		Taking the expectation over $D$, we find that the overall probability of recovering the root correctly is at least $1 - 2q \E{e^{-\frac{\delta^2q^2}{8} \lambda_2^2D}} $. Using the fact that the moment generating function of the Poisson distribution is $\E{e^{tD}} = e^{d(e^t-1)}$, and that for all feasible values of $\delta$, say $\delta < \frac{2}{q}$, we have that $e^{-\frac{\delta^2q^2}{8}} < 1 - \frac{\delta^2q^2}{32}$, we find that the probability is bounded below by 
		\[ 1 - 2q \E{e^{-\frac{\delta^2q^2}{8} \lambda_2^2D}} = 1 - 2q e^{d(e^{-\frac{\delta^2q^2}{8}\lambda_2^2}-1)} > 1 - 2q e^{-\frac{\delta^2q^2}{32} \lambda_2^2d}. \]
		As before, we have that the exact same argument provides the result for $\tilde E_k$. 
	\end{proof}
	
	Finally, as a corollary, we obtain an improved version of \cref{lemma: numChildCancel} using this new estimate. 
	\begin{corollary}\label{cor: numChildCancel}
		Suppose that $\sigma_\rho = i$. There exists a $C = C(\pi, q, \delta)$ such that if the communities are $\delta$-distinguishable by $\{\xi_i: \lambda_i^2d > C\}$, then for layer $L$ deep enough and any child $v$, 
		\[ \abs{\E{\Prob{\sigma_v = j \middle\vert \sigma_L^{(i)}}} - P_{ij}} \leq  4\sqrt2 q^2 \max \pi_k^{1/2} \max\pi_k^{1/2} \cdot \lambda_2 \cdot e^{-\frac{\delta^2q^2}{32} \lambda_2^2d}. \]
	\end{corollary}
	\begin{proof}
		The proof is the same as \cref{lemma: numChildCancel}, where we use \cref{prop: doubleMajority} instead of \cref{prop: majority} to bound $\Prob{\sigma_v = j \middle\vert \sigma_L^{(j)}}$. 
	\end{proof}
	
	\section{The Recursive Formula and its Derivatives}\label{section: recursion}
	In this section, we derive the recursive formula for $X_\rho^{(m)}(i)$ in terms of the vectors $X_v^{(m-1)}$ on its children. Then, we compute the derivative of this formula to aid in our further analysis. Finally, we prove the main proposition that will allow us to estimate the sizes of these derivatives in expectation, by analyzing the ratios of the main terms that appear in the recursive and derivative formula. 
	\subsection{Derivation from Bayes Rule}
	The formula is analogous to the standard tree recursion derived from Bayes' rule, but due to the involvement of matrices we include the proof for completeness.
	\begin{proposition}
		\[ X_\rho^{(m)} (i)  = \frac{ \pi(i) \prod_{j=1}^d \left[ PD_\pi^{-1} X_j^{(m-1)} \right]_i}{\sum_{\ell=1}^q  \pi(\ell) \prod_{j=1}^d \left[ PD_\pi^{-1} X_j^{(m-1)} \right]_\ell}. \]
	\end{proposition}
	\begin{proof}
		Starting from the definition and applying Bayes' Rule, we get
		\begin{align*}
			X_\rho^{(m)} (i) &= \Prob{\sigma_\rho = i \middle\vert \sigma_{L_m(\rho)} = L}. = \frac{\Prob{\sigma_{L_m(\rho)} = L \middle\vert \sigma_\rho = i}\Prob{\sigma_\rho = 	i}}{\sum_{\ell = 1}^q \Prob{\sigma_{L_m(\rho)} = L \middle\vert \sigma_\rho = \ell}\Prob{\sigma_\rho = \ell}}.
		\end{align*}
		From here, we look at the numerator, as the denominator is of a very similar form. By conditioning on the state of the states of the children, the distinct subtrees are independent and we can expand our probability into a product.
		\begin{align*}
			\Prob{\sigma_{L_m(\rho)} = L \middle\vert \sigma_\rho = i}\Prob{\sigma_\rho = i} &= \pi(i) 	\prod_{j=1}^d \sum_{k=1}^q \Prob{\sigma_{L_{m-1}(j)} = L\vert_{j} \middle\vert \sigma_{j} = k}\Prob{\sigma_{j} = k \middle\vert \sigma_\rho = i}. \\
			\intertext{Applying Bayes' Rule again to recover the vector $X_j^{(m-1)}$ defined on the children,}
			&= \pi(i) \prod_{j=1}^d \sum_{k=1}^q X_j^{(m-1)}(k)\cdot \Prob{\sigma_{L_{m-1}(j)} = L\vert_{j}} 	\cdot \frac{\Prob{\sigma_j = k \middle\vert \sigma_\rho = i}}{\Prob{\sigma_j = k}}  \\
			&= \pi(i) \prod_{j=1}^d \sum_{k=1}^q P_{ik} X_j^{(m-1)}(k) \cdot \frac{\Prob{\sigma_{L_{m-1}(j)} = 	L\vert_j}}{\pi(k)}. \\
			\intertext{Writing the inner sum as a particular entry of a matrix product, we get}
			&= \pi(i) \prod_{j=1}^d \Prob{\sigma_{L_{m-1}(j)} = L\vert_j}\left[ PD_\pi^{-1}X_j^{(m-1)}\right]_i. \\
			\intertext{Combining the product over the subtrees again by independence, we are left with }
			&= \Prob{\sigma_{L_m(\rho)} = L} \pi(i) \prod_{j=1}^d \left[ PD_\pi^{-1} X_j^{(m-1)} \right]_i.
		\end{align*}
		Returning to the original expression, we get that 
		\begin{align*}
			X_\rho^{(m)} (i)  &= \frac{\Prob{\sigma_{L_m(\rho)} = L} \pi(i) \prod_{j=1}^d \left[ PD_\pi^{-1} X_j^{(m-1)} \right]_i}{\sum_{\ell=1}^q \Prob{\sigma_{L_m(\rho)} = L} \pi(\ell) \prod_{j=1}^d \left[ PD_\pi^{-1} X_j^{(m-1)} \right]_\ell} = \frac{ \pi(i) \prod_{j=1}^d \left[ PD_\pi^{-1} X_j^{(m-1)} \right]_i}{\sum_{\ell=1}^q  \pi(\ell) 	\prod_{j=1}^d \left[ PD_\pi^{-1} X_j^{(m-1)} \right]_\ell}.
		\end{align*}
	\end{proof}

	From the above proposition, we define the recursive function in terms of the matrix $P$. 
	\begin{definition}\label{def: recursivefn}
		Let $f: \mathbb{R}^{d \times q} \rightarrow \mathbb{R}^q$ be defined by 
		\[ f_i(x_1, x_2, \ldots, x_d) = \frac{\pi(i)\prod_{j=1}^d (PD_\pi^{-1}X_j)_i}{\sum_{k=1}^q \pi(k)\prod_{j=1}^d (PD_{\pi}^{-1}X_j)_k}  \]
		where $f = (f_1, f_2, \ldots, f_q) $. 
	\end{definition}
	In order to make the derivative easier to analyze, we take the derivative with respect to a transformed version of the input vectors $x_i$. We first make a linear scaling and consider the vector $Y \coloneqq D_\pi^{-1}X$ and note that they satisfy the recursive expression 
	\[ Y_\rho^{(m)}(i) = \pi(i)^{-1} X_{\rho}^{(m)}(i) = \frac{\prod_{j=1}^d (PY_j^{(m-1)})_i}{\sum_{k=1}^q \pi(k) \prod_{j=1}^d (PY_j^{(m-1)})_k}. \]
	Noticing the similar structure in the numerator and denominator of this formula, we make the following definition.
	\begin{definition}\label{def: N}
		Let $N_i \coloneqq \prod_{j=1}^d (PY_j)_i$ and $N_{i, -t} \coloneqq \prod_{j \neq t} (PY_j)_i$.
	\end{definition}
	We can then rewrite \cref{def: recursivefn} as 
	\[ f_i = \frac{\pi_iN_i}{\sum_{k=1}^q \pi_kN_k}. \]
	We are also interested in working in the basis of eigenvectors of $P$. This will allow us to better isolate the actions of each eigenvalue. In particular, recall from \cref{fact: diagonalizing P} that we wrote $P = V^{-1} \Lambda V$ where $V = UD_p^{1/2}$ where $U$ is an orthogonal matrix. Along this line, we consider the transformed vectors $VY_\rho^{(m)}$ instead. 
	
	
	To analyze the derivative of the $f_i$, we first compute the partial derivative of each $N_i$ with respect to the individual entry $(VY_t)(\ell)$ of the transformed vector. Recall from \cref{def: N} that 
	\[ N_i = \prod_{j=1}^d (PY_j)_i = \prod_{j=1}^d (V^{-1}\Lambda VY_j)_i \]
	and note that we can express $V^{-1} = D_p^{-1}V^T$. This allows us to calculate that 
	\begin{align*}
		\frac{\partial N_i}{\partial (VY_t)(\ell)} &= \lambda_\ell (V^{-1})_{i\ell}\prod_{j \neq t} (PY_j)_i = \lambda_\ell \pi_i^{-1} V_{\ell i} N_{i, -t}.
	\end{align*}
	
	Finally, we can derive our expression for the derivative of the recursion with respect to $VY_\rho$. 
	\begin{align*}
		\frac{\partial (X_\rho^{(m)})_i}{\partial (VY_t^{(m-1)})_\ell} &= \frac{\left(\sum_{k=1}^q \pi_kN_k\right)\left(\lambda_\ell V_{\ell i}N_{i, -t}\right) - \left(\pi_iN_i\right)\left(\sum_{k=1}^q \lambda_\ell V_{\ell k}N_{k, -t}\right)}{\left(\sum_{k=1}^q \pi_kN_k\right)^2} \\
		&= \frac{\sum_{k= 1}^q \left[\pi_kN_k\lambda_\ell V_{\ell i}N_{i, -t} - \pi_iN_i\lambda_\ell V_{\ell k}N_{k, -t} \right]}{\left(\sum_{k=1}^q \pi_kN_k\right)^2} \\
		&= \lambda_\ell \cdot \frac{\sum_{k=1}^q \left(\pi_kV_{\ell i}N_kN_{i, -t} - \pi_iV_{\ell k}N_iN_{k, -t}\right)}{\left(\sum_{k=1}^q \pi_kN_k\right)^2}. 
	\end{align*}
	Most importantly, notice that in the numerator the terms of the form $N_iN_{i, -t}$ are canceled out by the coefficient. Thus, we only have products $N_iN_j$ for distinct indices $i \neq j$. On the other hand, in the denominator, we have the terms $N_i^2$. This is what we will take advantage of to show that the derivative is small with high probability. Thus, we derive our expression for the derivative as 
	\label{equatoin: derivative}
	\begin{equation}
		\frac{\partial (X_\rho^{(m)})_i}{\partial (VY_t^{(m-1)})_\ell} = \lambda_\ell \cdot \frac{\sum_{k=1, k \neq i}^q \left(\pi_kV_{\ell i}N_kN_{i, -t} - \pi_iV_{\ell k}N_iN_{k, -t}\right)}{\sum_{m,n=1}^q \pi_m\pi_nN_mN_n}.
	\end{equation}
	
	\begin{remark}
		When writing the functions $f_i$ or the partial derivatives $\frac{\partial f_i}{\partial (VY_t)_\ell}$ we will often omit the arguments for simplicity of notation. It should be understood that the arguments are typically some combination of the $X_v^{(n-1)}$ and $\tilde X_v^{(n-1)}$, which behave essentially equivalently under our analysis. In the event that specific inputs are designated they will be written out in the first instance. 
	\end{remark}
	
	\subsection{Analysis of the Derivative}
	As noted above, we use the fact that the squared terms $N_i^2$ show up in the denominator but not the numerator. To take advantage of this, we will show that typically, there is one $N_i$ that is much larger than the rest. If this is the case, then the squared term in the denominator will dominate the entire derivative, and thus we will obtain that in the typical case the derivative is very small. 
	\begin{proposition}\label{prop: Nanalysis}
		Suppose that $\frac{2\sqrt{2} \max \pi_k^{3/2} \max \pi_k^{-3/2}}{\xi^3} \lambda_2 < \frac{\delta^2q^2}{8}$ and the conditions of \cref{cor: numChildCancel} hold, and assume that the root has label 1. Let $d$ be the number of children of the root. Then for $m$ large enough, with probability at least $1-\exp{-C\lambda_2^2d}$, 
		\[ \frac{N_i}{N_1} \leq e^{-C\lambda_2^2 d}. \]
		An analogous inequality holds when the root label is replaced by any $i \in [q]$ and any vectors $X^{(m)}$ in the definition of $N$ are replaced by their noisy counterparts $\tilde X^{(m)}$. 
	\end{proposition}
	\begin{proof}
		The proof of this will be broken down into two steps. 
		\begin{enumerate}
			\item The ratio between the values of $N_1$ and $N_i$ with an exactly proportional distribution of children and perfect information on these children is large.
			\item The values of $N_i$ we from reconstructing based on leaves is is close to the ideal values from Step 1. 
		\end{enumerate}
		The second step will be further decomposed into a few steps, analyzing the different degree terms of a Taylor expansion. 
		
		\noindent \textit{Step 1.} 
		We first show a large ratio of approximate values of the $N_i$ and $N_1$. In an ideal world, we would have exactly $dP_{1j}$ children of community $j$, and on each of these children we would have perfect confidence in our guess. Under these scenario, our corresponding values of $N_i$ and $N_1$ are
		\[ A_{i1} = \prod_{j=1}^q  P_{ij}^{dP_{1j}} \text{ and } A_{11} = \prod_{j=1}^q P_{1j}^{dP_{1j}} \]
		respectively. Taking logs, we see that 
		\begin{align*}
			\log(A_{11}) - \log (A_{i1}) &= \sum_{j=1}^q dP_{1j}\log P_{1j} - \sum_{j=1}^q dP_{1j}\log P_{ij} = d\sum_{j=1}^q P_{1j} \log (\frac{P_{1j}}{P_{ij}}).
		\end{align*}
		This quantity is precisely the relative entropy from the distribution $P_i$ to the distribution $P_1$.  By assumption \cref{condition: rows far enough apart}, we know that these two distributions differ in $L^1$ norm by at least $\delta q \lambda_2$. Applying Pinsker's Inequality, we get that 
		\[ \sum_{j=1}^q P_{1j} \log (\frac{P_{1j}}{P_{ij}}) \geq \frac{1}{2} \cdot \norm{P_1 - P_j}^2 \geq \frac{\delta^2q^2}{2} \lambda_2^2. \]
		Returning to our original expressions, 
		\begin{align*}
			\frac{A_{11}}{A_{i1}} &= \exp(\log (A_{11}) - \log (A_{i1})) = \exp(d\sum_{j=1}^q P_{1j} \log (\frac{P_{1j}}{P_{ij}})) \geq \exp{\frac{\delta^2q^2}{2} \lambda_2^2 d}. 
		\end{align*}
		This completes step 1 of the proof. 
		
		\noindent \textit{Step 2.} In step 2, we would like to show that $N_i$ is close to $A_{i1}$ and similarly $N_1$ is close to $A_{11}$, so that an inequality of the same form continues to hold. We consider the ratio of $N_i$ to $A_{i1}$ by similarly taking logs. 
		\begin{align*}
			\log N_i - \log (A_{i1}) &= \sum_{k=1}^d \log((PY_k)_i) - \sum_{j=1}^q dP_{1j} \log (\frac{P_{ij}}{\pi_j}).
			\intertext{Suppose that there are exactly $n_j$ of the children in community $j$} 
			&= \sum_{k=1}^d \log((PY_k)_i) - \sum_{k=1}^d \log (\frac{P_{i\sigma(k)}}{\pi_{\sigma(k)}}) + \sum_{j=1}^q (n_j-dP_{1j})\log(\frac{P_{ij}}{\pi_j}).
		\end{align*}
		We analyze the difference of the first two terms together in the first part, then bound the third term in the second part of this step. 
		
		\noindent \textit{Step 2, Part 1.} The difference between the first two sums that we are interested in can then be written by Taylor's Theorem as 
		\begin{align*}
			\sum_{k=1}^d \log((PY_k)_i) - \sum_{k=1}^d \log (\frac{P_{i\sigma(k)}}{\pi_{\sigma(k)}})  &= \sum_{k=1}^d \left[\left( 1 - (PY_k)_i \right) - \left( 1 - (PY_k)_i \right)^2 + \frac{1}{t_k^3}\left( 1 - (PY_k)_i \right)^3 \right] \\
			&\quad - \sum_{k=1}^d \left[\left( 1 - \frac{P_{i\sigma(k)}}{\pi_{\sigma(k)}}\right) - \left( 1 - \frac{P_{i\sigma(k)}}{\pi_{\sigma(k)}}\right)^2 + \frac{1}{u_k^3}\left( 1 - \frac{P_{i\sigma(k)}}{\pi_{\sigma(k)}}\right)^3 \right] \\
			&= \sum_{k=1}^d \left( \frac{P_{i\sigma(k)}}{\pi_{\sigma(k)}} - (PY_k)_i\right) \\
			&\qquad - \sum_{k=1}^d \left[ \left( 1 - (PY_k)_i \right)^2 - \left( 1 - \frac{P_{i \sigma(k)}}{\pi_{\sigma(k)}}\right)^2 \right] \\
			&\qquad + \sum_{k=1}^d \left[ \frac{1}{t_k^3}\left( 1 - (PY_k)_i \right)^3 - \frac{1}{u_k^3}\left( 1 - \frac{P_{i \sigma(k)}}{\pi_{\sigma(k)}}\right)^3 \right] 
		\end{align*}
		where the $t_k$ are on the interval 1 to $(PY_k)_i$ and $u_k$ are on the interval 1 to $\frac{P_{i \sigma(k)}}{\pi_{\sigma(k)}}$. Note here that both these intervals are bounded below by $\frac{\xi}{\min \pi_k}$, and so we have a uniform upper bound on these random coefficients.
		
		\noindent \textit{Step 2, Part 1.1.}
		To analyze the linear difference term, we write
		\begin{align*}
			\epsilon_k \coloneqq (PY_k)_i - \frac{P_{i\sigma(k)}}{\pi_{\sigma(k)}} &= \sum_{j=1}^q \frac{P_{ij}}{\pi_j} \cdot \Prob{\sigma_k = j \middle\vert \sigma_L} - \sum_{j=1}^q \frac{P_{ij}}{\pi_j} \cdot \mathbbm{1}(\sigma_k = j) \\
			&= \sum_{j=1}^q \frac{P_{ij}}{\pi_j} \cdot (\Prob{\sigma_k = j \middle\vert \sigma_L} - \mathbbm{1}(\sigma_k = j)).
		\end{align*}
		Since the terms $\Prob{\sigma_k = j \middle\vert \sigma_L} - \mathbbm{1}(\sigma_k = j)$ sum to 0, we may subtract 1 from each of their coefficients without changing the value.
		\[ \epsilon_k = \sum_{j=1}^q \left(\frac{P_{ij}}{\pi_j}-1\right) \cdot (\Prob{\sigma_k = j \middle\vert \sigma_L} - \mathbbm{1}(\sigma_k = j)). \]
		In expectation, we have that 
		\begin{align*}
			\E{\epsilon_k} &= \E{\sum_{j=1}^q \left(\frac{P_{ij}}{\pi_j}-1\right) \cdot (\Prob{\sigma_k = j \middle\vert \sigma_L} - \mathbbm{1}(\sigma_k = j))} = \sum_{j=1}^q \left(\frac{P_{ij}}{\pi_j}-1\right) \cdot \E{\Prob{\sigma_k = j \middle\vert \sigma_L} - \mathbbm{1}(\sigma_k = j)}. 
		\end{align*}
		Recall that we are working conditioned on $\sigma_\rho = 1$. By \cref{cor: numChildCancel}, we know that the term \\ $\E{\Prob{\sigma_k = j \middle\vert \sigma_L} - \mathbbm{1}(\sigma_k = j)}$ is of the order $\lambda_2e^{-C\lambda_2^2d}$. Thus, our expectation is controlled in size by 
		\begin{align*}
			\abs{\E{\epsilon_k}} &\leq \sum_{j=1}^q \abs{\frac{P_{ij}}{\pi_j}-1} \cdot \abs{\E{\Prob{\sigma_k = j \middle\vert \sigma_L} - \mathbbm{1}(\sigma_k = j)}} \leq C\lambda_2^2 e^{-\frac{\delta^2q^2}{32}\lambda_2^2d}.
		\end{align*}
		At this point, we choose $\lambda_2^2d$ large enough so that this absolute value is at most $\lambda_2^2\frac{\delta^2q^2}{64}$. 
		Moreover, due to the factor of $\frac{P_{ij}}{\pi_j} - 1$ appearing in each term of the sum, and since the remaining factor $\Prob{\sigma_k = j \middle\vert \sigma_L} - \mathbbm{1}(\sigma_k = j)$ is bounded in absolute value by 1, we have that the residual term 
		\[ \abs{\epsilon_k} \leq \sqrt2 q\max  \pi_k^{1/2} \max \pi_k^{-1/2}\lambda_2.  \]
		Thus, when we have an independent sum of these residual terms over the children of the root, we may apply Hoeffding's Inequality and obtain a small error with high probability due to cancellation. To bound the sum in a high probability situation, we use the following form of Hoeffding's Inequality. 
		We have that 
		\begin{align*}
			\Prob{\abs{\sum_{k=1}^d \epsilon_k - \E{\epsilon_k}} \geq d\lambda_2^2 \cdot \frac{\delta^2q^2}{64}} &\leq 2\exp(-\frac{2d\lambda_2^4\delta^4q^4}{4096(\sqrt2 q\max  \pi_k^{1/2} \max \pi_k^{-1/2} \lambda_2)^2} )  \\
			&= 2\exp(-\frac{\delta^4q^2}{4096 \max  \pi_k \max \pi_k^{-1}}\cdot \lambda_2^2d).
		\end{align*}
		This moreover implies that we have 
		\[ \Prob{\abs{\sum_{k=1}^d \epsilon_k} \geq d\lambda_2^2 \cdot \frac{\delta^2q^2}{32} } \leq 2\exp(-\frac{\delta^4q^2}{4096 \max  \pi_k \max \pi_k^{-1}}\cdot \lambda_2^2d) \]
		due to our assumption on the expectation $\E{\epsilon_k}$. This concludes Step 2, Part 1.1. 
		
		\noindent \textit{Step 2, Part 1.2.} For the second order terms from the expansion, we perform a very similar analysis. Following the same ideas, we first estimate the expectation of 
		\begin{align*}
			\epsilon_k^{(2)} &= (1 - (PY_k)_i)^2 - \left( 1 - \frac{P_{i \sigma(k)}}{\pi_{\sigma(k)}} \right)^2. 
		\end{align*}
		Since the terms $\Prob{\sigma_k = j \middle\vert \sigma_L}$ and $\mathbbm{1}(\sigma_k = j)$ sum to 1, we can pull the $-1$ inside the summation that defines $(PY_k)_i$ and $\frac{P_{i \sigma(k)}}{\pi_{\sigma(k)}}$ and attach them to the coefficients. This gives
		\begin{align*}
			((PY_k)_i-1)^2 - \left(\frac{P_{i\sigma(k)}}{\pi_{\sigma(k)}}-1\right)^2 &= \left(\sum_{j=1}^q \left(\frac{P_{ij}}{\pi_j}-1\right) \cdot \Prob{\sigma_k = j \middle\vert \sigma_L}\right)^2 - \left(\sum_{j=1}^q \left(\frac{P_{ij}}{\pi_j} - 1\right) \cdot \mathbbm{1}(\sigma_k = j)\right)^2. \\
			\intertext{Squaring out, we get}
			&= \sum_{j, j'=1}^q \left(\frac{P_{ij}}{\pi_j}-1\right) \cdot \left(\frac{P_{ij'}}{\pi_{j'}}-1\right) \\
			&\qquad \qquad  \cdot (\Prob{\sigma_k = j \middle\vert \sigma_L}\Prob{\sigma_k = j' \middle\vert \sigma_L} - \mathbbm{1}(\sigma_k = j)\mathbbm{1}(\sigma_k = j')). \\
			\intertext{Extracting the terms for which $j = j'$, }
			&= \sum_{j=1}^q \left(\frac{P_{ij}}{\pi_j}-1\right)^2 \cdot (\Prob{\sigma_k = j \middle\vert \sigma_L}^2 - \mathbbm{1}(\sigma_k = j)) \\
			&\qquad + \sum_{j, j'=1}^q \left(\frac{P_{ij}}{\pi_j}-1\right) \cdot \left(\frac{P_{ij'}}{\pi_{j'}}-1\right) \cdot (\Prob{\sigma_k = j \middle\vert \sigma_L}\Prob{\sigma_k = j' \middle\vert \sigma_L}). 
		\end{align*}
		Each of the terms $\left(\frac{P_{ij}}{\pi_j}-1\right)$ holds a factor of $\lambda_2$, so our entire sum will have a contribution of $\lambda_2^2$ from these terms. Note that this implies we have a uniform bound on $\epsilon_k^{(2)}$ of the form $C\lambda_2^2$, since the probabilities and indicators are of magnitude at most 1. As for a more precise estimate of the probabilities, we consider then conditioned on $\sigma_k$. We know that 
		\begin{align*}
			0 \geq \E{\Prob{\sigma_k = j \middle\vert \sigma_L}^2 - \mathbbm{1}(\sigma_k = j) \middle\vert \sigma_k = j} &= \E{\Prob{\sigma_k=j \middle\vert \sigma_L}^2 \middle\vert \sigma_k = j} - 1 \\
			&\geq \E{\Prob{\sigma_k=j \middle\vert \sigma_L} \middle\vert \sigma_k = j}^2 - 1 \\
			&\geq (1 - Ce^{-C\lambda_2^2d})^2 - 1 \\
			&\geq -Ce^{-C\lambda_2^2d}
		\end{align*}
		and for $j' \neq j$,
		\begin{align*}
			0 \leq \E{\Prob{\sigma_k = j \middle\vert \sigma_L}^2 - \mathbbm{1}(\sigma_k = j) \middle\vert \sigma_k = j'} &= \E{\Prob{\sigma_k=j \middle\vert \sigma_L}^2 \middle\vert \sigma_k = j'} \\
			&\leq \E{\Prob{\sigma_k=j \middle\vert \sigma_L} \middle\vert \sigma_k = j}  \\
			&\leq Ce^{-C\lambda_2^2d}.
		\end{align*}
		Similarly, for 
		\begin{align*}
			0 \leq \E{\Prob{\sigma_k = j \middle\vert \sigma_L}\Prob{\sigma_k = j' \middle\vert \sigma_L} \middle\vert \sigma_k = j} &\leq \E{\Prob{\sigma_k = j' \middle\vert \sigma_L} \middle\vert \sigma_k = j} \leq Ce^{-C\lambda_2^2d}
		\end{align*}
		and for $\ell \neq j$, 
		\begin{align*}
			0 \leq \E{\Prob{\sigma_k = j \middle\vert \sigma_L}\Prob{\sigma_k = j' \middle\vert \sigma_L} \middle\vert \sigma_k = \ell} &\leq \E{\Prob{\sigma_k = j \middle\vert \sigma_L} \middle\vert \sigma_k = \ell} \leq Ce^{-C\lambda_2^2d}.
		\end{align*}
		Combining all of these estimates, we obtain from the form above that 
		\[ \abs{\E{\epsilon_k^{(2)}}} \leq C\lambda_2^2 e^{-C\lambda_2^2d}. \]
		This matches with what we obtained in Step 2, Part 1.1. Once again, at this point we ensure that $\lambda_2^2d$ is large enough so that this quantity is at most $\lambda_2^2\frac{\delta^2q^2}{64}$. With this assumption, we may apply Hoeffding's Inequality to show that
		\begin{align*}
			\Prob{\abs{\sum_{k=1}^d \epsilon_k^{(2)} - \E{\epsilon_k^{(2)}}} \geq d\lambda_2^2 \cdot \frac{\delta^2q^2}{64}} &\leq 2\exp(-\frac{2d\lambda_2^4\delta^4q^4}{4096(C\lambda_2^2)^2} ) = 2\exp(-Cd).
		\end{align*}
		This moreover implies that we have 
		\[ \Prob{\abs{\sum_{k=1}^d \epsilon_k^{(2)}} \geq d\lambda_2^2 \cdot \frac{\delta^2q^2}{32} } \leq 2\exp(-Cd) \]
		due to our assumption on the expectation $\E{\epsilon_k^{(2)}}$. This concludes Step 2, Part 1.2. 
		
		\noindent \textit{Step 2, Part 1.3}. For the third order terms from the expansion, we simply apply a uniform bound. Recall from \cref{lemma: Pclosetopi} that 
		\[ \abs{P_{i\sigma(k)} - \pi_{\sigma(k)}} \leq \sqrt{2} \max \pi_k^{1/2} \max \pi_k^{-1/2} \lambda_2. \]
		Diving both sides by $\pi_{\sigma(k)}$, we obtain the upper bound 
		\[ \abs{\frac{P_{i\sigma(k)} }{ \pi_{\sigma(k)}} - 1}  \leq \sqrt{2} \max \pi_k^{1/2} \max \pi_k^{-1/2} \cdot \frac{1}{\min \pi_k} \lambda_2. \]
		Thus, cubing and taking the sum over all children, and taking the uniform upper bound on $\frac{1}{u_k^2}$, we get 
		\[ \abs{\sum_{k=1}^d \frac{1}{u_k^3} \left( 1 - \frac{P_{i \sigma(k)}}{\pi_{\sigma(k)}} \right)^3}  \leq \frac{2\sqrt{2} \max \pi_k^{3/2} \max \pi_k^{-3/2}}{\xi^3} \lambda_2^3d.  \]
		In particular, this third order term is of size $O(\lambda_2^3d)$, which is what will suffice. 
		
		Similarly, we obtain a bound on the term $(PY_k)_i - 1$ as $(PY_k)_i$ is simply a convex combination of the $\frac{P_{i\sigma(k)}}{\pi_{\sigma(k)}}$. Expanding as above, we can write 
		\begin{align*}
			\abs{(PY_k)_i - 1} &= \abs{\sum_{j=1}^q \left(\frac{P_{ij}}{\pi_j} \cdot \Prob{\sigma_k = j \middle\vert \sigma_L}\right) - 1}. \\
			\intertext{Pulling the $-1$ into the summation as explained previously, }
			&= \abs{\sum_{j=1}^q \left(\frac{P_{ij}}{\pi_j} - 1 \right) \cdot \Prob{\sigma_k = j \middle\vert \sigma_L}}. \\
			\intertext{By the triangle inequality, }
			&\leq \sum_{j=1}^q \abs{\frac{P_{ij}}{\pi_j} - 1} \cdot \Prob{\sigma_k = j \middle\vert \sigma_L} \\
			&\leq  \sum_{j=1}^q \left( \sqrt{2} \max \pi_k^{1/2} \max \pi_k^{-1/2} \cdot \frac{1}{\min \pi_k} \lambda_2 \right) \cdot \Prob{\sigma_k = j \middle\vert \sigma_L} \\
			&= \sqrt{2} \max \pi_k^{1/2} \max \pi_k^{-1/2} \cdot \frac{1}{\min \pi_k} \lambda_2.
		\end{align*}
		Thus, by the same calculation as above, we can obtain the same uniform bound 
		\[ \abs{\sum_{k=1}^d \frac{1}{t_k^3} \left( 1 - (PY_k)_i \right)^3}  \leq \frac{2\sqrt{2} \max \pi_k^{3/2} \max \pi_k^{-3/2}}{\xi^3} \lambda_2^3d. \]
		This concludes Step 2, Part 1.3.
		
		Combining these intermediate steps to estimate of the first two terms in our expansion and complete Step 2, Part 1, we obtain 
		\begin{align*}
			\abs{\sum_{k=1}^d \log((PY_k)_i) - \sum_{k=1}^d \log (\frac{P_{i\sigma(k)}}{\pi_{\sigma(k)}})} 
			&\leq \abs {\sum_{k=1}^d \left( \frac{P_{i\sigma(k)}}{\pi_{\sigma(k)}} - (PY_k)_i\right)} \\
			&\qquad + \abs{ \sum_{k=1}^d \left[ \left( 1 - (PY_k)_i \right)^2 - \left( 1 - \frac{P_{i \sigma(k)}}{\pi_{\sigma(k)}}\right)^2 \right] } \\
			&\qquad + \abs{\sum_{k=1}^d \left[ \frac{1}{t_k^3}\left( 1 - (PY_k)_i \right)^3 - \frac{1}{u_k^3}\left( 1 - \frac{P_{i \sigma(k)}}{\pi_{\sigma(k)}}\right)^3 \right] } \\
			&\leq d\lambda_2^2 \cdot \frac{\delta^2q^2}{32} + d\lambda_2^2 \frac{\delta^2q^2}{32} + \frac{2\sqrt{2} \max \pi_k^{3/2} \max \pi_k^{-3/2}}{\xi^3} \lambda_2^3d \\
			&= \left( \frac{\delta^2q^2}{16} + \frac{2\sqrt{2} \max \pi_k^{3/2} \max \pi_k^{-3/2}}{\xi^3} \lambda_2 \right) \lambda_2^2d.
		\end{align*}
		
		\noindent \textit{Step 2, Part 2}. To conclude Step 2, we need to control the remaining term in the sum 
		\[ \sum_{j=1}^q (n_j-dP_{1j})\log(\frac{P_{ij}}{\pi_j}). \]
		Applying Taylor's Theorem to the log, we find that 
		\begin{align*}
			\log(\frac{P_{ij}}{\pi_j}) &= \log(1 + (\frac{P_{ij}}{\pi_j} - 1)) = \frac{1}{t}(\frac{P_{ij}}{\pi_j} - 1)
		\end{align*}
		for some $t$ between 1 and $\frac{P_{ij}}{\pi_j}$. Note here that the lower bound of this interval is again at least $\frac{\xi}{\min \pi_k}$ so in particular the term $t^{-1}$ is bounded above by a constant. Evaluating with an absolute value, by Hoeffding's Inequality, we see that with probability at least $1-2qe^{-2\alpha^2\lambda_2^2d}$,
		\begin{align*}
			\abs{\sum_{j=1}^q (n_j-dP_{1j})\log(\frac{P_{ij}}{\pi_j})} &\leq \sum_{j=1}^q \abs{n_j-dP_{1j}} \cdot \abs{\frac{1}{t_j}(\frac{P_{ij}}{\pi_j}-1)}  \\
			&\leq \sum_{j=1}^q \alpha\lambda_2d \cdot \frac{\xi}{\min \pi_k} \cdot \sqrt2 \max \pi_k^{1/2} \max\pi_k^{-1/2}\lambda_2 \\
			&\leq \alpha \cdot \frac{q\min \pi_k\cdot \sqrt2 \max \pi_k^{1/2} \max\pi_k^{-1/2}}{\xi} \cdot  \lambda_2^2d.
		\end{align*}
		All together, we have shown that with high probability, i.e. probability at least $1 - Ce^{-C\lambda_2^2d}$, 
		\begin{align*}
			\abs{\log N_i - \log A} &\leq \abs{\sum_{k=1}^d \log((PY_k)_i) - \sum_{k=1}^d \log (\frac{P_{i\sigma(k)}}{\pi_{\sigma(k)}}) + \sum_{j=1}^q (n_j-dP_{ij})\log(\frac{P_{ij}}{\pi_j})} \\
			&\leq  \abs{\sum_{k=1}^d \log((PY_k)_i) - \sum_{k=1}^d \log (\frac{P_{i\sigma(k)}}{\pi_{\sigma(k)}})} + \abs{\sum_{j=1}^q (n_j-dP_{ij})\log(\frac{P_{ij}}{\pi_j})} \\
			&\leq \left(\alpha \frac{q\min \pi_k\cdot \sqrt2 \max \pi_k^{1/2} \max\pi_k^{-1/2}}{\xi} +  \frac{\delta^2q^2}{16} + \frac{2\sqrt{2} \max \pi_k^{3/2} \max \pi_k^{-3/2}}{\xi^3} \lambda_2  \right) \cdot  \lambda_2^2d.
		\end{align*}
		Thus, taking the exponential, we find the desired bound that 
		\[ e^{-C^*\lambda_2^2d} \leq \frac{N_i}{A_{i1}} \leq e^{C^*\lambda_2^2d} \]
		where in this expression $C^*$ is the coefficient obtained in the above estimate. This concludes Step 2.
		
		Finally, we may estimate the ratio of the two terms we set out to check: $N_1$ and $N_i$. We have that 
		\begin{align*}
			\frac{N_1}{N_i} &= \frac{N_1}{A_{11}} \cdot \frac{A_{11}}{A_{i1}} \cdot \frac{A_{i1}}{N_i} \geq e^{-C\lambda_2^2d} \cdot e^{\frac{\delta^2q^2}{2}\lambda_2^2d} \cdot e^{-C\lambda_2^2d} = e^{\frac{\delta^2q^2}{8} \cdot \lambda_2^2d}.
		\end{align*}
		Here, we choose $\alpha$ so that the exponent works out as we expect it to. 
		\begin{align*}
			\frac{\delta^2q^2}{2} - 2\left(\alpha \frac{q\min \pi_k\cdot \sqrt2 \max \pi_k^{1/2} \max\pi_k^{-1/2}}{\xi} +  \frac{\delta^2q^2}{16} + \frac{2\sqrt{2} \max \pi_k^{3/2} \max \pi_k^{-3/2}}{\xi^3} \lambda_2  \right) \geq \frac{\delta^2q^2}{8} 
		\end{align*}
		As long as $\frac{2\sqrt{2} \max \pi_k^{3/2} \max \pi_k^{-3/2}}{\xi^3} \lambda_2 < \frac{\delta^2q^2}{8}$ we can pick $\alpha$ small enough in terms of $q, \pi, \xi$ such that the given inequality holds. This provides the desired high probability bound. 
	\end{proof}
	
	\begin{remark}\label{remark: weaken constraint}
		Regarding the additional constraint of $\frac{2\sqrt{2} \max \pi_k^{3/2} \max \pi_k^{-3/2}}{\xi^3} \lambda_2 < \frac{\delta^2q^2}{8}$, we remark that it is weak enough for there to still be an interesting class of  matrices satisfying these conditions. However, in the event that this constraint is not satisfactory, it is easy to check  that with the same argument as above, by taking the Taylor expansion to the $K$th term and bounding the difference of each power analogously, it is possible to achieve a constraint of the form $C\lambda_2^{K-2} < \frac{\delta^2q^2}{8}$. This may be useful in some scenarios, but we do not investigate this further. 
	\end{remark}
	
	\begin{corollary}
		\cref{prop: Nanalysis} holds even when we vary the input vector $X_t$ for any one fixed $t$ arbitrarily. 
	\end{corollary}
	\begin{proof}
		There are a couple of ways to see this corollary. The first is to note that changing any single vector $X_t$ cannot change the corresponding $N_i$ by much. In particular, it can only multiply or divide the term by a constant factor. 
		
		An alternative approach which may be clearer is to note the effect on the application of Hoeffding's inequality. In the case that we do not have control over $X_t$, we simply apply Hoeffding's inequality on the remaining $d-1$ children. Then, re-introducing child $t$ can only add or subtract a bounded amount. Thus, our difference in logarithms changes by at most a constant, which corresponds to the above multiplication or division by a constant when the exponential is taken. 
		
		In either case, we are allowed to arbitrarily change one input vector $X_t$, which will be essential in completing the analysis in the next chapter. 
	\end{proof}
	
	To close out the section, we apply the above \cref{prop: Nanalysis} to obtain a bound on the expected first partial derivative of our recursive functions $f_i$. This will set the example for the analysis of all other derivatives we encounter, at which point some details may be omitted for the sake of clarity. 
	
	\begin{proposition}\label{prop: derivative example}
		Suppose that $\frac{2\sqrt{2} \max \pi_k^{3/2} \max \pi_k^{-3/2}}{\xi^3} \lambda_2 < \frac{\delta^2q^2}{8}$ and the conditions of \cref{cor: numChildCancel} hold. Suppose the root has $d$ children. We then have that for $m$ deep enough, 
		\[ \abs{\E{ \frac{\partial f}{\partial (VY_t)_\ell}(X_1, \ldots, X_d)  \middle\vert \tilde X_t}} \leq C\lambda_2^2 e^{-C\lambda_2^2 d}. \]
	\end{proposition}
	\begin{proof}
		We would like to bound the derivative independent of the input vector for child $t$. Since the term $(PY_j)_i$ is bounded below by $\frac{\xi}{\max \pi_k}$ for all $i, j$, we can factor out the $PY_t$ factor from each $N_i$. Recall from \cref{equatoin: derivative} that we expressed the derivative in the form 
		\begin{align*}
			\frac{\partial (X_\rho^{(m)})_i}{\partial (VY_t^{(m-1)})_\ell} &= \lambda_\ell \cdot \frac{\sum_{k=1, k \neq i}^q \left(\pi_kV_{\ell i}N_kN_{i, -t} - \pi_iV_{\ell k}N_iN_{k, -t}\right)}{\sum_{m,n=1}^q \pi_m\pi_nN_mN_n} \\
			&\leq \lambda_\ell \cdot \frac{2q \max \pi_k N_{1} \max_{j \neq 1} N_{j}}{\xi \pi_1^2 N_{1}^2} = \frac{2q \max \pi_k}{\xi \pi_1^2} \cdot \lambda_2 \cdot \left( \max_{j \neq 1} \frac{N_{j}}{N_{1}}\right).
		\end{align*}
		Note, we have written the bound in the form when conditioned on the true community of the root being a 1. This is because in this case, we will have that $N_1 >> N_j$ for all other $j$. This will be true for any conditioning on the value of $\sigma_\rho$, and so the bound we obtain in this conditional setting will extend to the general unconditioned setting as well. By \cref{prop: Nanalysis}, we know that with probability at least $1-(q-1)\exp(-C(\delta, \xi, q, \pi) \lambda_2^2 d)$ the term in the right-most parentheses is at most $\exp(-\frac{\delta^2q^2}{8} \cdot \lambda_2^2d)$. In the remaining low probability situation, we can make a  global bound on the derivative. 
		\begin{align*}
			\lambda_\ell \cdot \frac{\sum_{k=1, k \neq i}^q \left(\pi_kV_{\ell i}N_kN_{i, -t} - \pi_iV_{\ell k}N_iN_{k, -t}\right)}{\sum_{m,n=1}^q \pi_m\pi_nN_mN_n} &\leq \lambda_\ell \cdot \frac{2\max \pi_k \sum_{m , n=1}^q N_mN_{n}}{\xi \min \pi_k^2 \sum_{m,n=1}^q N_mN_n} = \frac{2\max \pi_k}{\xi \min \pi_k^2} \cdot \lambda_\ell = \lambda_\ell \cdot \frac{2\max \pi_k}{\xi \min\pi_k^2}.
		\end{align*}
		All together, we have that 
		\begin{align*}
			\abs{\E{ \frac{\partial f}{\partial (VY_t)_\ell} \middle\vert \tilde X_t, \sigma_\rho=1}} &\leq \frac{2q \max \pi_k}{\xi \min \pi_k^2} \lambda_2 \exp(-\frac{\delta^2q^2}{8} \cdot \lambda_2^2d) + (q-1)\exp(-C \lambda_2^2 d)\lambda_\ell \cdot \frac{2\max \pi_k}{\xi \min\pi_k^2} \\
			&\leq \lambda_2C_1\exp(-C_2\lambda_2^2d)
		\end{align*}
		where $C_1$ and $C_2$ depend only on $\delta, \xi, q$, and $\pi$. As mentioned before, since an analogous argument holds for the cases when $\sigma_\rho$ is in any other community, we have overall that 
		\[ \abs{\E{ \frac{\partial f}{\partial (VY_t)_\ell}  \middle\vert \tilde X_t}} \leq \lambda_2C_1 \exp(-C_2\lambda_2^2d). \]
	\end{proof}
	
	\subsection{Other Derivatives}\label{section: higher derivatives}
	In this section, we make a few remarks about higher order derivatives of $f_i$ and derivatives of functions related to $f_i$ that will contribute to our analysis. First, on writing out the form of the higher order derivatives, we can see that each partial derivative introduces an addition factor of $\lambda$. Thus, we will have that a $k$th order partial derivative of $f_i$ will be on the order of $\lambda_2^k$. We now have to consider the main term of the derivative. Once again, just as with the first partial derivatives, we will have one of the $N_i$'s dominating the remaining $N_j$ with high probability. We also have that the numerator is a lower degree polynomial in each of the $N_i$ than the denominator. In total, this explanation implies that for $\arg\max N_i = m$, 
	\[ \partial^\alpha f_i = O\left(\lambda_2^{\abs{\alpha}} \cdot \frac{\sum_{i\neq m} N_i}{N_m}\right). \]
	Analogously to the analysis of the first partial derivatives above, we have that $\abs{\partial^\alpha f_i} \leq \lambda_2^{\abs{\alpha}} e^{-C \lambda_2^2d}$ with high probability and $\abs{\partial^\alpha f_i} \leq C\lambda_2^{\abs{\alpha}}$ uniformly. 
	
	We will also be interested in derivatives of functions such as $f_if_j$ and $f_i^2$. Note that by the product rule, we have 
	\[ \partial f_if_j = f_i\partial f_j + f_j\partial f_i \quad \text{ and } \quad \partial f_i^2 = 2f_i \partial f_i.  \]
	Since our functions $f_i$ recursively define a probability, we will identically have that $f_i \in [0, 1]$. Thus, in these partial derivatives, we obtain the similar bounds 
	\[ \abs{\partial f_if_j} = \abs{\partial f_j} + \abs{\partial f_i} \quad \text{ and } \quad \abs{\partial f_i^2} = 2 \abs{ \partial f_i }. \]
	Applying the bounds we obtained in the last section, we find that the partial derivatives of $f_if_j$ and $f_i^2$ obey the same estimates up to constant factors. This will also be important to our final analysis. 
	
	\section{Contraction of Noisy and Non-Noisy Probabilities}\label{section: contraction}
	Denote $(\mathcal{E}_n)_{ij} = \E{X_i^{(n)}(j) - \tilde X_i^{(n)}(j) \middle \vert \sigma_\rho = i}$ to be the matrix of errors between the noisy and non-noisy estimates. Along this line, let $\epsilon_n = \max_{ij} \abs{(\mathcal{E}_n)_{ij}}$ be the maximum absolute error on the $n$th sublevel. The goal of this chapter will be to show that this discrepancy $\epsilon_n$ tends to 0 as $n$ tends to infinity. Before proceeding to the calculation, we present a series of technical lemmas that will be essential to the analysis. 
	
	\subsection{Technical Lemmas}
	
	To start off, we present the consequence of cancellation that will appear often in the computation. This cancellation is the key to introducing factors of $\lambda_2$ where they may not appear trivially. For the first example of this, we can compute that $\E{X_i(j) - \tilde X_i(j) \middle\vert \sigma_\rho = k} $ is on the order of $\lambda_2 \epsilon_n$. In particular, when conditioned on the label of a parent, the error term gains an extra factor of $\lambda_2$. 
	\begin{lemma}\label{lemma: difference}
		For any $i$ a child of $\rho$ and $j, k \in [q]$, $\abs{\E{X_i(j) - \tilde X_i(j) \middle\vert \sigma_\rho = k}} \leq C \lambda_2 \epsilon_n. $
	\end{lemma}
	
	\begin{proof}
		\begin{align*}
			\abs{\E{X_i(j) - \tilde X_i(j)\middle\vert \sigma_\rho = k}} &= \abs{\sum_{\ell = 1}^q \E{X_i(j) - \tilde X_i(j) \middle\vert \sigma_i = \ell} \Prob{\sigma_i = \ell \middle\vert \sigma_\rho = k}} = \abs{\sum_{\ell = 1}^q (E_n)_{\ell j} \left(P_{k \ell} - \pi_\ell \right)} \\
			&\leq \sum_{\ell = 1}^q \abs{(E_n)_{\ell j} \left(P_{k \ell} - \pi_\ell \right)} \leq  \sum_{k \neq j, k=1}^q \epsilon_n \cdot \sqrt{2} q \max \pi_k^{1/2} \max \pi_k^{-1/2} \lambda_2 \\
			&= \sqrt{2} q^2 \max \pi_k^{1/2} \max \pi_k^{-1/2}\lambda_2 \epsilon_n.
		\end{align*}
		Here, in the second equality we can subtract $\pi_k$ from each probability since the coefficients are conditional probabilities that will average to 0 over all possibilities. The following calculation verifies this:
		\begin{align*}
			\sum_{k=1}^q \E{X_i(j) - \tilde X_i(j) \middle\vert \sigma_i = k} \pi_k &= \E{X_i(j) - \tilde X_i(j)} = \pi_j - \pi_j = 0.
		\end{align*}
	\end{proof}
	
	The second lemma provides a crucial bound on the second moment of the difference. It relies on the nature of $X$ and $\tilde X$ defined as conditional probabilities with increasing levels of information, a relationship which will be explored further in the next section. 
	
	\begin{lemma}\label{lemma: differenceSquared}
		$ \E{(X(i) - \tilde X(i))^2} = \pi_i \cdot (\mathcal{E}_n)_{ii}. $
	\end{lemma}
	\begin{proof}
		We first write out $X$ and $\tilde X$ conditioning on the noisy leaves on the $n$th sublevel. 
		\begin{align*}
			\tilde X(i) &= \Prob{\sigma_\rho = i \middle\vert \tau(n)} = \sum_{\tilde z} \Prob{\sigma_\rho = i \middle\vert \tau = \tilde z}  \Prob{\tau = \tilde z } \\
			X(1) &= \Prob{\sigma_\rho = i \middle\vert \sigma(n)} = \Prob{\sigma_\rho = i \middle\vert \sigma(n), \tau(n)} \\
			&= \sum_{\tilde z} \Prob{\sigma_\rho = i \middle\vert \sigma(n), \tau(n) = \tilde z} \Prob{\tau(n) = \tilde z }.
		\end{align*}
		As these two sums are taken over the same index set, we can evaluate the differences between each of the corresponding terms. In particular, set 
		\[ W = \Prob{\sigma_\rho = i \middle\vert \tau(n) = \tilde z, \sigma(n)}. \]
		This is a random variable, which has expectation 
		\[ \E{W} = \Prob{\sigma_\rho = i \middle\vert \tau(n) = \tilde z}. \]
		Conveniently, these are precisely the terms that appear in the sums for $X$ and $\tilde X$ respectively. Thus, we turn out attention to deriving an equivalent expression for $\E{W - \E{W} \middle\vert \sigma_\rho = i}$ which provides a relation to \\ $\E{(W - \E{W})^2}$. We evaluate $\E{W \middle\vert \sigma_\rho = i}$ using a similar technique to Lemma 2.2 in \cite{S:11}. 
		
		First, we expand the expectation over all possibilities of $\sigma(n)$, the $n$th sublevel. 
		\[ \E{W \middle\vert \sigma_\rho = i} = \sum_z \Prob{\sigma_\rho = i \middle\vert \sigma(n) = z, \tau = \tilde z} \Prob{\sigma(n) = z \middle\vert \sigma_\rho = i, \tau = \tilde z}. \]
		Using Bayes' Rule, we can rewrite the second conditional probability within the summation.
		\[ \E{W \middle\vert \sigma_\rho = i} = \sum_z \Prob{\sigma_\rho = i \middle\vert \sigma(n) = z, \tau(n) = \tilde z} \cdot \frac{\Prob{\sigma_\rho = i \middle\vert \sigma(n) = z, \tau(n) = \tilde z}\Prob{\sigma(n) = z \middle\vert \tau(n) = \tilde z}}{\Prob{\sigma_\rho = i \middle\vert \tau = \tilde z}}. \]
		Grouping like terms together, we get the expression
		\[ \E{W \middle\vert \sigma_\rho = i} = \sum_z  \Prob{\sigma_\rho = i \middle\vert \sigma(n) = z, \tau(n) = \tilde z}^2 \cdot \frac{\Prob{\sigma(n) = z \middle\vert \tau(n) = \tilde z}}{\Prob{\sigma_\rho = i \middle\vert \tau(n) = \tilde z}}. \]
		Rewriting the summation in terms of an expectation over $\sigma(n)$ again, 
		\[ \E{W \middle\vert \sigma_\rho = i} = \E{\Prob{\sigma_\rho = i \middle\vert \sigma(n), \tau(n) = \tilde z}^2} \cdot \frac{1}{\E{W}}. \]
		Writing all terms in terms of $W$ and manipulating the expression to contain $\Var{W}$, we get the relation
		\[ \E{W \middle\vert \sigma_\rho = i} = \frac{\E{W^2}}{\E{W}} = \frac{\E{W}^2 + \Var{W}}{\E{W}}  = \E{W} + \frac{\Var{W}}{\E{W}}. \]
		Subtracting the expectation to the other side, since it is a constant we get that
		\[ \E{W - \E{W}\middle\vert \sigma_\rho = i} = \frac{\E{(W-\E{W})^2}}{\E{W}}. \]
		This is exactly the form of expression that we desired. Returning to our original expression, we get that 
		\begin{align*}
			\E{X(i) - \tilde X(i) \middle\vert \sigma_\rho = i} &= \sum_{\tilde z} \E{W - \E{W} \middle\vert \sigma_\rho = i} \Prob{\tau(n) = \tilde z \middle\vert \sigma_\rho = i} \\
			&= \sum_{\tilde z} \E{(W - \E{W})^2} \cdot \frac{\Prob{\tau(n) = \tilde z \middle\vert \sigma_\rho  =i}}{\Prob{\sigma_\rho = i \middle\vert \tau(n) = \tilde z}} \\
			&= \sum_{\tilde z} \E{(W - \E{W})^2} \cdot \frac{\Prob{\tau(n) = \tilde z}}{\Prob{\sigma_\rho = i}} \\
			&= \frac{1}{\Prob{\sigma_\rho = i}} \cdot \E{(X(i) - \tilde X(i))^2}.
		\end{align*}
		Using our notation, we get that 
		\[ \E{(X(i) - \tilde X(i))^2} = \Prob{\sigma_\rho = i} \cdot (\mathcal{E}_n)_{ii}. \]
	\end{proof}
	From the above, we obtain the immediate bound in terms of $\epsilon_n$. 
	\begin{corollary}
		$ \E{(X(i) - \tilde X(i))^2} \leq \pi_i \epsilon_n. $
	\end{corollary}
	
	In similar spirit, we can obtain an analogous result for the other second order terms as well. 
	\begin{lemma}\label{lemma: covariance}
		$ \E{(X(i) - \tilde X(i))(X(j) - \tilde X(j))} = \pi_j (\mathcal{E}_n)_{ij}. $
	\end{lemma}
	\begin{proof}
		We perform the an analogous calculation to \cref{lemma: differenceSquared}, which slightly different conditioning. Here, we set 
		\[ W_i = \Prob{\sigma_\rho = i \middle\vert \tau(n) = \tilde z, \sigma(n)} \]
		with 
		\[ \E{W_i} = \Prob{\sigma_\rho = i \middle\vert \tau(n) = \tilde z}. \]
		In this case, we analyze $\E{W_i - \E{W_i} \middle\vert \sigma_\rho = j}$, which provides a relation to 
		\[ \E{(W_i - \E{W_i})(W_j - \E{W_j})}. \] The steps as taken above proceed as follows. 
		\begin{align*}
			\E{W_i \middle\vert \sigma_\rho = j} &= \sum_z \Prob{\sigma_\rho = i \middle\vert \sigma(n) = z, \tau = \tilde z} \Prob{\sigma(n) = z \middle\vert \sigma_\rho = j, \tau = \tilde z} \\
			&= \sum_z \Prob{\sigma_\rho = i \middle\vert \sigma(n) = z, \tau(n) = \tilde z}  \cdot \frac{\Prob{\sigma_\rho = j \middle\vert \sigma(n) = z, \tau(n) = \tilde z}\Prob{\sigma(n) = z \middle\vert \tau(n) = \tilde z}}{\Prob{\sigma_\rho = j \middle\vert \tau = \tilde z}} \\
			&= \sum_z  \Prob{\sigma_\rho = i \middle\vert \sigma(n) = z, \tau(n) = \tilde z}\Prob{\sigma_\rho = j \middle\vert \sigma(n) = z, \tau(n) = \tilde z} \cdot \frac{\Prob{\sigma(n) = z \middle\vert \tau(n) = \tilde z}}{\Prob{\sigma_\rho = j \middle\vert \tau(n) = \tilde z}} \\
			&= \E{\Prob{\sigma_\rho = i \middle\vert \sigma(n), \tau(n) = \tilde z}\Prob{\sigma_\rho = j \middle\vert \sigma(n), \tau(n) = \tilde z}} \cdot \frac{1}{\E{W_j}} \\
			\E{W_i \middle\vert \sigma_\rho = j} &= \frac{\E{W_iW_j}}{\E{W_j}} = \frac{\E{W_i}\E{W_j} + \Cov{W_i}{W_j}}{\E{W_j}}  = \E{W_i} + \frac{\Cov{W_i}{W_j}}{\E{W_j}} \\
			\E{W_i - \E{W_i}\middle\vert \sigma_\rho = j} &= \frac{\E{(W_i-\E{W_i})(W_j - \E{W_j})}}{\E{W_j}}.
		\end{align*}
		Returning to our original expression, we get that 
		\begin{align*}
			\E{X(i) - \tilde X(i) \middle\vert \sigma_\rho = j} &= \sum_{\tilde z} \E{W_i - \E{W_i} \middle\vert \sigma_\rho = j} \Prob{\tau(n) = \tilde z \middle\vert \sigma_\rho = j} \\
			&= \sum_{\tilde z} \E{(W_i-\E{W_i})(W_j - \E{W_j})} \cdot \frac{\Prob{\tau(n) = \tilde z \middle\vert \sigma_\rho  =j}}{\Prob{\sigma_\rho = j \middle\vert \tau(n) = \tilde z}} \\
			&= \sum_{\tilde z} \E{(W_i-\E{W_i})(W_j - \E{W_j})} \cdot \frac{\Prob{\tau(n) = \tilde z}}{\Prob{\sigma_\rho = j}} \\
			&= \frac{1}{\Prob{\sigma_\rho = j}} \cdot \E{(W_i-\E{W_i})(W_j - \E{W_j})}
		\end{align*}
		which implies that
		\[ \E{(X(i) - \tilde X(i))(X(j) - \tilde X(j))} = \Prob{\sigma_\rho = j} \cdot (\mathcal{E}_n)_{ij}. \]
	\end{proof}
	Once again, we obtain an immediate bound in absolute value in terms of $\epsilon_n$ as a consequence. 
	\begin{corollary}
		$ \abs{\E{(X(i) - \tilde X(i))(X(j) - \tilde X(j))}} \leq \pi_j \epsilon_n. $
	\end{corollary}
	
	As an application to demonstrate the utility of the above lemmas, we compute the expectation of a quantity that will naturally appear in our calculations. As seen from the following lemma, we obtain a bound on the squared difference of our transformed vectors of the same nature. 
	\begin{lemma}\label{lemma: differenceSquaredVY}
		$ \E{(VY(i) - V\tilde Y(i))^2} \leq C\epsilon_n. $
	\end{lemma}
	\begin{proof}
		We can expand the difference  $VY(i) - V\tilde Y(i)$ as 
		\[ VY(i) - V\tilde Y(i) = \sum_{j=1}^q V_{ij}(Y(j) - \tilde Y(j)). \]
		Squaring, we get the expansion 
		\begin{align*}
			(VY(i) - V\tilde Y(i))^2 &= \sum_{j=1}^q \sum_{k=1}^q V_{ij}V_{ik} (Y(j) - \tilde Y(j))(Y(k) - \tilde Y(k)) \\
			&= \sum_{j=1}^q V_{ij}^2(Y(j) - \tilde Y(j))^2 + \sum_{j \neq k} V_{ij}V_{ik}(Y(j) - \tilde Y(j))(Y(k) - \tilde Y(k)).
		\end{align*}
		Since the entries of $Y$ are at most $\max \pi^{-1}$ times the entries of $X$, the above \cref{lemma: differenceSquared} shows that the squared terms in the first sum are of order at most $C\epsilon_n$ in expectation, where $C$ depends only on $\pi$. Moreover, the coefficients $V_{ij}$ are entries of vectors of norm $\sqrt{\pi_j}$, and so are also bounded in terms of $\pi$. Similarly, we may bound the terms in the second sum in expectation by \cref{lemma: covariance}, and so we have the desired bound. 
	\end{proof}
	
	\subsection{Martingale Interpretation}
	
	We introduce an interpretation of the random variables $X$ and $\tilde X$ as a martingale that will motivate and simplify our calculations. Consider the three-term sequence $M_0  = \tilde X_k(i), M_1 = X_k(i), M_2 = \mathbbm{1}(\sigma_k = i)$. We check that this sequence defines a martingale. Indeed, we can rewrite them from definition as $M_0 = \Prob{\sigma_k = i \middle\vert \tau_{L_n}}, M_1 = \Prob{\sigma_k = i \middle\vert \tau_{L_n}, \sigma_{L_n}}, M_2 = \Prob{\sigma_k = i \middle\vert \tau_{L_n}, \sigma_{L_n}, \sigma_k}$ and note that these correspond to increasing levels of information in the conditioning, forming an exposure martingale with respect to $\tau_{L_n}, \sigma_{L_n}, \sigma_k$. 
	
	The first lemma from this interpretation allows us to simplify a squared difference into a difference of squares. This will be useful as it will reduce the computation of a double sum to the simpler computation of a single sum. 
	\begin{lemma}
		$ \E{(X_\rho(i) - \tilde X_\rho(i))^2} = \E{X_\rho(i)^2 - \tilde X_\rho(i)^2}. $
	\end{lemma}
	\begin{proof}
		\begin{align*}
			\E{(X_k(i) - \tilde X_k(i))^2} &= \E{X_k(i)^2  - 2 X_k(i)\tilde X_k(i) + \tilde X_k(i)^2} \\
			&= \E{X_k(i)^2} - 2 \E{\E{X_k(i) \tilde X_k(i) \middle\vert \tau_{L_n}}} + \E{\tilde X_k(i)^2} \\
			&= \E{X_k(i)^2} - 2 \E{\tilde X_k(i)\E{X_k(i)  \middle\vert \tau_{L_n}}} + \E{\tilde X_k(i)^2} \\
			&= \E{X_k(i)^2} - 2 \E{\tilde X_k(i)^2} + \E{\tilde X_k(i)^2} \\
			&=  \E{X_k(i)^2 - \tilde X_k(i)^2}.
		\end{align*}
	\end{proof}
	
	Similarly, we can simplify the covariance terms as well, in particular the remaining order two terms of the difference between $X$ and $\tilde X$. The proof is very similar. 
	\begin{lemma}\label{lemma: martingale simplify 2}
		$ \E{(X_\rho(i) - \tilde X_\rho(i))(X_\rho(j)) - \tilde X_\rho(j)} = \E{X_\rho(i)X_\rho(j) - \tilde X_\rho(i)\tilde X_\rho(j)}. $
	\end{lemma}
	
	The next lemma informs us when certain expectations that appear in our calculation are identically positive or identically negative. This will be helpful in our bounds when the absolute value of certain terms come into play. 
	\begin{lemma}\label{lemma: martingale positive}
		$\E{(X_k(i) - \tilde X_k(i))\mathbbm{1}(\sigma_k = i) \middle\vert  \tau_{L_n}} \geq 0 $ and 
		$ \E{(X_k(i) - \tilde X_k(i))\mathbbm{1}(\sigma_k \neq i) \middle\vert  \tau_{L_n}} \leq 0. $
	\end{lemma}
	\begin{proof}
		The quantity we want to analyze is $\E{(M_1-M_0)M_2 \middle\vert \tau_{L_n}}$. 
		\begin{align*}
			\E{(M_1-M_0)M_2 \middle\vert \tau_{L_n}} &= \E{\E{(M_1-M_0)M_2 \middle\vert \tau_{L_n}, \sigma_{L_n}}\middle\vert \tau_{L_n}} \\
			&= \E{(M_1-M_0)\E{M_2 \middle\vert \tau_{L_n}, \sigma_{L_n}} \middle\vert \tau_{L_n}} \\
			&= \E{(M_1-M_0)M_1 \middle\vert \tau_{L_n}} \\
			&= \E{(M_1-M_0)M_1 \middle\vert \tau_{L_n}} - \E{(M_1-M_0)M_0 \middle\vert \tau_{L_n}} \\
			&= \E{(M_1-M_0)^2 \middle\vert \tau_{L_n}} \\
			&\geq 0.
		\end{align*}
	\end{proof}
	
	\subsection{Contraction}
	We can finally synthesize all of our work and prove the contraction of the errors. We would like to relate $\epsilon_{n+1}$ to $\epsilon_n$ by estimating $(\mathcal{E}_{n+1})_{ij} = \E{X_\rho^{(n+1)}(i) - \tilde X_\rho^{(n+1)}(i) \middle\vert \sigma_\rho = j}$ in relation to the entries of $\mathcal{E}_n$. Here, we use our technical lemmas to adjust the quantity we analyze in order to aid the calculation. By \cref{lemma: covariance}, we have that 
	\[ (\mathcal{E}_{n+1})_{ij} = \pi_j^{-1} \E{(X_\rho^{(n+1)}(i) - \tilde X_\rho^{(n+1)}(i)) (X_\rho^{(n+1)}(j) - \tilde X_\rho^{(n+1)}(j))}. \]
	Applying \cref{lemma: martingale simplify 2}, we obtain that 
	\[ (\mathcal{E}_{n+1})_{ij} = \pi_j^{-1} \E{X_\rho^{(n+1)}(i)X_\rho^{(n+1)}(j) - \tilde X_\rho^{(n+1)}(i) \tilde X_\rho^{(n+1)}(j)} \]
	From our recursive formula, we can express 
	\[  X_\rho^{(n+1)}(i) = f_i(X_1^{(n)}, \ldots, X_D^{(n)}) \text{ and } \tilde X_\rho^{(n+1)}(i) = f_i(\tilde X_1^{(n)}, \ldots, \tilde X_D^{(n)}). \]
	We denote the input vector of all children $X^{(n)}$ and its noisy counterpart $\tilde X^{(n)}$. Also, denote $X^{(n)}_{-k}$ to be the input vector of all children except for the $k$th child. Thus, considering the function $f_if_j$, we can express 
	\[ (\mathcal{E}_{n+1})_{ij} = \pi_j^{-1}\E{f_if_j(X^{(n)}) - f_if_j(\tilde X^{(n)})}. \]
	We turn out attention to this difference of $f_if_j$ evaluated at exact versus noisy leaves. To analyze this difference, we expand by telescoping to isolate the changes on each of the individual children, then Taylor expand. By telescoping, we obtain
	\[ f_if_j(X^{(n)}) - f_if_j(\tilde X^{(n)}) = \sum_{k=1}^{D} f_if_j(Z_{k-1}^{(n)}) - f_if_j(Z_{k}^{(n)}) \]
	where $Z_k^{(n)}$ is the input vector $X_i^{(n)}$ for $1 \leq i \leq D-k$ and $\tilde X_i^{(n)}$ for $D-k+1 \leq i \leq D$. Similarly denote $(Z_k^{(n)})_{-k}$ to be this vector without the input from the $k$th child.  In particular, the input vector where we have changed to noisy leaves on the last $k$ children. Applying the expectation once again, and using the Tower Rule to condition on the number of children $D$, we get 
	\begin{align*}
		(\mathcal{E}_{n+1})_{ij} &= \pi_j^{-1}\E{f_if_j(X^{(n)}) - f_if_j(\tilde X^{(n)})} \\
		&= \pi_j^{-1}\E{\sum_{k=1}^{D} f_if_j(Z_{k-1}^{(n)}) - f_if_j(Z_{k}^{(n)})} \\
		&= \pi_j^{-1}\E{\E{\sum_{k=1}^{D} f_if_j(Z_{k-1}^{(n)}) - f_if_j(Z_{k}^{(n)}) \middle\vert D}} \\
		&= \pi_j^{-1}\E{ \sum_{k=1}^{D} \E{ f_if_j(Z_{k-1}^{(n)}) - f_if_j(Z_{k}^{(n)}) \middle\vert D}}.
	\end{align*}
	Thus, in order to bound our desired quantity $(\mathcal{E}_{n+1})_{ij}$ it suffices to estimate the quantities \[ \E{ f_if_j(Z_{k-1}^{(n)}) - f_if_j(Z_{k}^{(n)}) \middle\vert D}. \]
	
	With this form, we can then move to our next step of Taylor expanding the difference. Using a multivariate form of Taylor's Theorem, we can write 
	\begin{align*}
		f_if_j(Z_{k-1}^{(n)}) &= f_if_j(Z_{k}^{(n)}) + (VY_{k}^{(n)} - V\tilde Y_{k}^{(n)})^T \frac{\partial f_if_j}{\partial VY_{k}}(Z_k^{(n)}) + (VY_{k}^{(n)} - V\tilde Y_{k}^{(n)})^T H (VY_{k}^{(n)} - V\tilde Y_{k}^{(n)})
	\end{align*}
	where the partial derivatives are evaluated at $Z_{k}^{(n)}$ and the matrix $H$ is defined in terms of the second partial derivatives of $f_if_j$ with respect to the entries of $(VY_{k})_m$ on the interval between $X_{k}^{(n)}$ and $\tilde X_{k}^{(n)}$. In particular, 
	\[ H_{ab} = \int_0^1 (1-t)\frac{\partial^2 f_if_j}{\partial (VY_{k})_a \partial (VY_{k})_b}((Z_k^{(n)})_{-k}, \tilde X_{k} + t(X_{k} - \tilde X_{k})) \, dt. \]
	Note that 
	\[ \abs{H_{ab}} \leq \max_t \abs{\frac{\partial^2 f_if_j}{\partial (VY_{k})_a \partial (VY_{k})_b}((Z_k^{(n)})_{-k}, \tilde X_{k} + t(X_{k} - \tilde X_{k}))}. \]
	Since our bound on the derivative allows for independence from one of the input vectors corresponding to a specific child $k$, this quantity will be able to be bounded independent of $X_{k}$ and $\tilde X_{k}$. 
	
	We first handle the quadratic term, as it turns out to be the easier term to control. Notice that by the symmetry of partial derivatives, we have that $H$ is a symmetric matrix. This fact leads to the following lemma.
	\begin{lemma}\label{lemma: quadratic term}
		Let $\lambda_{max}(H)$ be the largest absolute eigenvalue of $H$. We have that 
		\[ \abs{(VY_{k}^{(n)} - V\tilde Y_{k}^{(n)})^T H (VY_{k}^{(n)} - V\tilde Y_{k}^{(n)})} \leq \norm{VY_{k}^{(n)} - V \tilde Y_{k}^{(n)}}_2^2 \lambda_{max}(H). \]
	\end{lemma}
	\begin{proof}
		By the above remark, we have that $H$ is a real symmetric matrix, so by the Spectral Theorem we can diagonalize $H = U^T \Lambda U$ where $U$ is a unitary matrix. Then, 
		\begin{align*}
			(VY_{k}^{(n)} - V\tilde Y_{k}^{(n)})^T H (VY_{k}^{(n)} - V\tilde Y_{k}^{(n)}) = (VY_{k}^{(n)} - V\tilde Y_{k}^{(n)})^T U^T \Lambda U (VY_{k}^{(n)} - V\tilde Y_{k}^{(n)}).
		\end{align*}
		Let $\vec{v} =  U (VY_{k}^{(n)} - V\tilde Y_{k}^{(n)})$ and note that since $U$ is unitary, $\norm{\vec{v}}_2 = \norm{VY_{k}^{(n)} - V\tilde Y_{k}^{(n)}}_2$. Making the substitution of notation, we have our expression as $\vec{v}^T \Lambda \vec{v}$. Notice that the diagonal matrix $\Lambda$ can contribute a factor of at most $\lambda_{max}(H)$ to the norm $\vec{v}^T\vec{v}$, and so we obtain that 
		\[ \abs{\vec{v}^T\Lambda \vec{v}} \leq \norm{\vec{v}}_2^2 \lambda_{max}(H) = \norm{VY_{k}^{(n)} - V\tilde Y_{k}^{(n)}}_2^2 \lambda_{max}(H). \]
	\end{proof}
	
	In expectation, this term can then be estimated by the following lemma. 
	\begin{lemma}\label{lemma: Taylor quadratic}
		Suppose that $\frac{2\sqrt{2} \max \pi_k^{3/2} \max \pi_k^{-3/2}}{\xi^3} \lambda_2 < \frac{\delta^2q^2}{8}$ and the conditions of \cref{cor: numChildCancel} hold. Then for $n$ sufficiently large
		\[ \E{\norm{VY_{k}^{(n)} - V \tilde Y_{k}^{(n)}}_2^2 \lambda_{max}(H)  \middle\vert D} \leq C\lambda_2^2 e^{-C\lambda_2^2D } \epsilon_n. \]
	\end{lemma}
	\begin{proof}
		We condition within the expectation on $\sigma_\rho$ and $\tau_{L_n(k)}$ so that the two terms are conditionally independent. We can then factor the inner terms and continue with the bound. 
		\begin{align*}
			\E{\norm{VY_{k}^{(n)} - V \tilde Y_{k}^{(n)}}_2^2 \lambda_{max}(H)  \middle\vert D} &= \E{\E{\norm{VY_{k}^{(n)} - V \tilde Y_{k}^{(n)}}_2^2 \lambda_{max}(H)  \middle\vert D, \sigma_\rho, \tau_{L_n(k)}}\middle\vert D} \\
			&= \E{\E{A \middle\vert D, \sigma_\rho, \tau_{L_n(k)}}\E{ B  \middle\vert D, \sigma_\rho, \tau_{L_n(k)}}\middle\vert D}
			\intertext{where $A = \norm{VY_{k}^{(n)} - V \tilde Y_{k}^{(n)}}_2^2$ and $B = \lambda_{max}(H)$. Notice that since $A$ is identically positive, we can upper bound this quantity by factoring our a global maximum of the expectation of $B$.}
			&\leq \norm{\E{B \middle\vert D, \sigma_\rho, \tau_{L_n(k)}}}_\infty \E{\E{A \middle\vert D, \sigma_\rho, \tau_{L_n(k)}} \middle\vert D} \\
			&= \norm{\E{B \middle\vert D, \sigma_\rho, \tau_{L_n(k)}}}_\infty \E{A \middle\vert D } 
			\intertext{Here, a random variable defined only on a child $k$ is independent from the number of children, and so we can remove the conditioning on $D$.}
			&= \norm{\E{B \middle\vert D, \sigma_\rho, \tau_{L_n(k)}}}_\infty \E{A} 
		\end{align*}
		
		We now bound these terms separately. Notice that 
		\[ A = \norm{VY_{k}^{(n)} - V \tilde Y_{k}^{(n)}}_2^2 = \sum_{m=1}^q (VY_{k}^{(n)}(m) - V \tilde Y_{k}^{(n)}(m))^2. \]
		By  \cref{lemma: differenceSquaredVY}, we have that each of these terms in the sum have expectation bounded by $C \epsilon_n$. Thus, overall, we have that $ \E{A} \leq q \cdot C \epsilon_n$.
		
		To control $\lambda_{max}(H)$, we note that the largest absolute eigenvalue of a real matrix is at most the maximum $L^1$ norm of a row of the matrix. For our matrix $H$, this is 
		\[ \sum_{b=1}^q \abs{H_{ab}} \leq q \max_t \abs{\frac{\partial^2 f_if_j}{\partial (VY_{k})_a \partial (VY_{k})_b}((Z_k^{(n)})_{-k}, \tilde X_{k} + t(X_{k} - \tilde X_{k}))}.  \]
		Now, following the idea of \cref{prop: derivative example} and the discussion in \cref{section: higher derivatives}, we note that the second partial derivative can be bounded with probability at least $1 - e^{-C\lambda_2^2D}$ independently of the input from child $k$ by $C\lambda_2^2e^{-C\lambda_2^2D}$ and $C\lambda_2^2$ uniformly. Thus, with high probability we have that 
		\[ \lambda_{max}(H) \leq q C\lambda_2^2e^{-C\lambda_2^2D} \] 
		and it is at most $q C \lambda_2^2$ uniformly. Then, in expectation
		\[ \norm{\E{\lambda_{max}(H) \middle\vert D, \sigma_\rho, \tau_{L_n(k)}}}_\infty \leq C\lambda_2^2e^{-C\lambda_2^2D}. \]
		All together, we get the desired result. 
	\end{proof}
	
	Now, we may turn to handling the linear term of the expansion, which will require a bit more care. First, we re-express the term as the sum of some nicer terms. Recall the form of the linear term 
	\begin{align*}
		(VY_{k}^{(n)} - V\tilde Y_{k}^{(n)})^T \frac{\partial f_if_j}{\partial VY_{k}}(Z_k^{(n)}) &= \sum_{m=1}^q (VY_{k}^{(n)}(m) - V \tilde Y_{k}^{(n)}(m)) \frac{\partial f_i f_j}{\partial VY_{k}(m)}(Z_k^{(n)}). \\
		\intertext{Expanding the definition of $VY$, we obtain a double sum of the form } 
		&= \sum_{m=1}^q \sum_{\ell=1}^q V_{m \ell} (Y_{k}^{(n)}(\ell) -  \tilde Y_{k}^{(n)}(\ell)) \frac{\partial f_i f_j}{\partial VY_{k}(m)}(Z_k^{(n)}) \\
		&= \sum_{m=1}^q \sum_{\ell=1}^q V_{m \ell}\pi_\ell^{-1} (X_{k}^{(n)}(\ell) -  \tilde X_{k}^{(n)}(\ell)) \frac{\partial f_i f_j}{\partial VY_{k}(m)}(Z_k^{(n)}).
	\end{align*}
	Since there are only $q^2$ terms in total, which is a constant, and the coefficients $V_{m\ell}$ are bounded in size in terms of $\pi$, it suffices to control the expectation of the summands $(X_{k}^{(n)}(\ell) -  \tilde X_{k}^{(n)}(\ell)) \frac{\partial f_i f_j}{\partial VY_{k}(m)} $. 
	
	At this point, we Taylor expand again to isolate the term in which we need to exploit the most cancellation. In this case, we Taylor expand around the $k$th input vector centered at the constant vector $\vec{c}$ which has all entries $\frac{1}{q}$. Using the multivariate Taylor's Theorem again on the function $g(v_{k}) = \frac{\partial f_if_j}{\partial VY_{k}(m)}((Z_k^{(n)})_{-k}, v_k)$, we can write 
	\begin{align*}
		g(\tilde X_{k}^{(n)}) &= g(\vec{c})  + (V \tilde Y_{k}^{(n)} - VD_\pi^{-1}\vec{c})^T \frac{\partial g}{\partial VY_{k}}(t\tilde X_k^{(n)} + (1-t)\vec{c}) \\
		&= g(\vec{c}) + \sum_{a=1}^q \left(V\tilde Y_{k}^{(n)}(a) - VD_\pi^{-1}\vec{c}(a)\right) \frac{\partial g}{\partial VY_{k}(a)}(t\tilde X_k^{(n)} + (1-t)\vec{c}) \\
		&= g(\vec{c}) + \sum_{a=1}^q\sum_{b=1}^q V_{ab}\left(\tilde Y_{k}^{(n)}(b) - D_\pi^{-1}\vec{c}(b)\right) \frac{\partial g}{\partial VY_{k}(a)}(t\tilde X_k^{(n)} + (1-t)\vec{c}) \\
		&= g(\vec{c}) + \sum_{a=1}^q\sum_{b=1}^q V_{ab}\pi_b^{-1} \left(\tilde X_{k}^{(n)}(b) - \frac{1}{q}\right) \frac{\partial g}{\partial VY_{k}(a)} (t\tilde X_k^{(n)} + (1-t)\vec{c}).
	\end{align*}
	Here and in the following computations, the first partial derivative of $g$ is evaluated at a point between $\vec{c}$ and $\tilde X_{k}^{(n)}$, so the argument may be omitted. This will not affect our bounds, as they are independent of this input vector, but it will define the dependence between certain random variables. Substituting this into our expression, we are interested in 
	\[ (X_{k}^{(n)}(\ell) - \tilde X_{k}^{(n)}(\ell)) g(\vec{c}) + \sum_{a=1}^q\sum_{b=1}^q V_{ab}\pi_b^{-1}  (X_{k}^{(n)}(\ell) - \tilde X_{k}^{(n)}(\ell)) \left(\tilde X_{k}^{(n)}(b) - \frac{1}{q}\right) \frac{\partial g}{\partial VY_{k}(a)}. \]
	As before, we now work to control the sizes of each of the summands. 
	\begin{lemma}\label{lemma: taylor taylor constant}
		Suppose that $\frac{2\sqrt{2} \max \pi_k^{3/2} \max \pi_k^{-3/2}}{\xi^3} \lambda_2 < \frac{\delta^2q^2}{8}$ and the conditions of \cref{cor: numChildCancel} hold. Then for $n$ sufficiently large
		\[ \abs{\E{(X_{k}^{(n)}(\ell) - \tilde X_{k}^{(n)}(\ell)) g(\vec{c}) \middle\vert D}} \leq C\lambda_2^2 e^{-C\lambda_2^2D} \epsilon_n.  \]
	\end{lemma}
	\begin{proof}
		Notice that the function $g$, as we are evaluating it at a constant point, depends only on the subtrees rooted at children other than $k$. Thus, by conditioning on the label of the root, we completely remove the dependence between the difference $(X_{k}^{(n)}(\ell) - \tilde X_{k}^{(n)}(\ell)) $ and $g(\vec{c})$. This allows us to factor the expectation with this additional conditioning. 
		\begin{align*}
			\abs{\E{(X_{k}^{(n)}(\ell) - \tilde X_{k}^{(n)}(\ell)) g(\vec{c}) \middle\vert D}} &= \abs{\E{\E{(X_{k}^{(n)}(\ell) - \tilde X_{k}^{(n)}(\ell)) g(\vec{c}) \middle\vert D, \sigma_\rho} \middle\vert D}} \\
			&= \abs{\E{\E{A \middle\vert D, \sigma_\rho} \E{B \middle\vert D, \sigma_\rho} \middle\vert D}}
			\intertext{where $A = (X_{k}^{(n)}(\ell) - \tilde X_{k}^{(n)}(\ell)) $ and $B = g(\vec{c})$.}
			&\leq \E{\abs{\E{A \middle\vert D, \sigma_\rho}} \abs{\E{B \middle\vert D, \sigma_\rho}} \middle\vert D}.
		\end{align*}
		Looking at the term $\E{A \middle\vert D, \sigma_\rho}$, we notice that again $A$ is independent from $D$ and so we simply have $\E{A \middle\vert \sigma_\rho}$. Then, by \cref{lemma: difference}, we have that this term is bounded in absolute value by $C\lambda_2 \epsilon_n$. 
		
		Now consider the term $\E{B \middle\vert D, \sigma_\rho}$. This is a first partial derivative that we have analyzed before. We know from \cref{prop: derivative example} that it is bounded by $C\lambda_2e^{-C\lambda_2^2D}$ in expectation. Recall that this is since we have the bound with high probability independent of the input vector for the $k$th child. Thus, combining these estimates, we find that 
		\[ \abs{\E{(X_{k}^{(n)}(\ell) - \tilde X_{k}^{(n)}(\ell)) g(\vec{c}) \middle\vert D}} \leq C\lambda_2^2 e^{-C\lambda_2^2D} \epsilon_n. \]
	\end{proof}
	
	Finally, we have the following estimate for the last type of summand. 
	\begin{lemma}\label{lemma: taylor linear linear}
		Suppose that $\frac{2\sqrt{2} \max \pi_k^{3/2} \max \pi_k^{-3/2}}{\xi^3} \lambda_2 < \frac{\delta^2q^2}{8}$ and the conditions of \cref{cor: numChildCancel} hold. Then for $n$ sufficiently large 
		\[ \E{(X_{k}^{(n)}(\ell) - \tilde X_{k}^{(n)}(\ell)) \left(\tilde X_{k}^{(n)}(b) - \frac{1}{q}\right) \frac{\partial g}{\partial VY_{k}(a)}  \middle\vert D} \leq C\lambda_2^2 e^{-C\lambda_2^2D} \epsilon_n. \]
	\end{lemma}
	\begin{proof}
		This proof will be very similar in spirit to the bounds we have already completed. The main difference in this case is that we have a dependence on $\tilde X_{k}^{(n)}$ in the partial derivative. Thus, in order to factor the expectation, we need to condition on $\tau_{L_n(k)}$ as well. Note that in this case, we will condition on $\sigma_{k}$ the label of child $k$ rather than the label of the root. Since this vertex also separates this particular subtree from the others, it also suffices for independence. 
		\begin{align*}
			\E{AB \middle\vert D} &= \E{\E{AB  \middle\vert D, \sigma_{k}, \tau_{L_n(k)}} \middle\vert D} = \E{\E{A  \middle\vert D, \sigma_{k}, \tau_{L_n(k)}}\E{B  \middle\vert D, \sigma_{k}, \tau_{L_n(k)}} \middle\vert D} 
		\end{align*}
		where $A = (X_{k}^{(n)}(\ell) - \tilde X_{k}^{(n)}(\ell)) \left(\tilde X_{k}^{(n)}(b) - \frac{1}{q}\right)$ and $B = \frac{\partial g}{\partial VY_{k}(a)} $. Taking absolute values, we obtain an expression as above in \cref{lemma: taylor taylor constant}. 
		
		In the case of $\E{B  \middle\vert D, \sigma_{k}, \tau_{L_n(k)}}$, we essentially have a second partial derivative of our initial function $f_if_j$. In this case we have conditioned on the state of a child instead of the state of the root. However, due to reversibility, and the fact that both $D$ and $\tau_{L_n(k)}$ are independent from the state of the root conditioned on $\sigma_{k}$, we have by the reversibility of the broadcast process that the root is distributed according to the row $P_{\sigma_{k}}$. Moreover, once conditioned on the state of the root, the dependence of all terms in the derivative on $\sigma_{k}$ disappears, and so we are simply in an average situation of the various labels of the root. Since our bound holds for all values of $\sigma_\rho$, it will also hold in this averaged situation. Thus, in particular, we have that this term is always bounded by $C\lambda_2^2 e^{-C\lambda_2^2D}$. 
		
		Now we move to the trickier case of $\E{A  \middle\vert D, \sigma_{k}, \tau_{L_n(k)}}$. Here, we have no dependence on $D$, so the conditioning will be removed. Referring to the proof of \cref{lemma: covariance} we have the equality that 
		\begin{align*}
			\E{A  \middle\vert \sigma_{k} = \ell', \tau_{L_n(k)}} &= \E{A(X_{k}^{(n)}(\ell') - \tilde X_{k}^{(n)}(\ell'))  \middle\vert  \tau_{L_n(k)}}
			\intertext{Substituting back in for $A$, we get }
			&= \E{(X_{k}^{(n)}(\ell) - \tilde X_{k}^{(n)}(\ell)) \left(\tilde X_{k}^{(n)}(b) - \frac{1}{q}\right)(X_{k}^{(n)}(\ell') - \tilde X_{k}^{(n)}(\ell'))  \middle\vert  \tau_{L_n(k)}}
			\intertext{Taking an absolute value, and noting that $\abs{\tilde X_{k}^{(n)}(b) - \frac{1}{q}}$, we can upper bound by the remaining terms }
			&\leq \E{\abs{(X_{k}^{(n)}(\ell) - \tilde X_{k}^{(n)}(\ell)) (X_{k}^{(n)}(\ell') - \tilde X_{k}^{(n)}(\ell'))}  \middle\vert  \tau_{L_n(k)}} 
			\intertext{Using the fact that $\frac{1}{2} (a^2 + b^2) \geq \abs{ab}$, we can upper bound this again by }
			&\leq \frac{1}{2} \E{(X_{k}^{(n)}(\ell) - \tilde X_{k}^{(n)}(\ell))^2 +  (X_{k}^{(n)}(\ell') - \tilde X_{k}^{(n)}(\ell'))^2  \middle\vert  \tau_{L_n(k)}}.
		\end{align*}
		Denote this by $\E{A' \middle\vert \tau_{L_n(k)}}$. Note importantly that this expectation is identically positive, which will be useful in concluding the computation. 
		\begin{align*}
			\abs{\E{AB \middle\vert D}} &\leq \E{\E{A' \middle\vert \tau_{L_n(k)}} \abs{\E{B  \middle\vert D, \sigma_{k}, \tau_{L_n(k)}}} \middle\vert D}
			\intertext{Now, since $\E{A' \middle\vert \tau_{L_n(k)}}$ is identically positive, we may upper bound this by factoring out the global maximum of the expectation of $B$}
			&\leq \norm{\E{B  \middle\vert D, \sigma_{k}, \tau_{L_n(k)}}}_\infty \E{\E{A' \middle\vert \tau_{L_n(k)}} \middle\vert D} \\
			&= \norm{\E{B  \middle\vert D, \sigma_{k}, \tau_{L_n(k)}}}_\infty \E{A'  \middle\vert D} 
			\intertext{Since $A'$ is independent from $D$, and including the bound on $B$ described above, we end up with }
			&\leq C\lambda_2^2 e^{-C\lambda_2^2D} \E{A'} 
			\intertext{Now, by \cref{lemma: differenceSquared}, we have that $\E{A'} \leq C\epsilon_n$ which allows us to conclude with the desired bound}
			&\leq C\lambda_2^2 e^{-C\lambda_2^2D} \epsilon_n.
		\end{align*}
	\end{proof}
	
	Combining the estimates of \cref{lemma: Taylor quadratic}, \cref{lemma: taylor taylor constant}, \cref{lemma: taylor linear linear}, we get the following estimate.
	\begin{lemma}
		Suppose that $\frac{2\sqrt{2} \max \pi_k^{3/2} \max \pi_k^{-3/2}}{\xi^3} \lambda_2 < \frac{\delta^2q^2}{8}$ and the conditions of \cref{cor: numChildCancel} hold. Then for $n$ sufficiently large
		\[ \abs{ \E{f_if_j(Z_{k-1}^{(n)}) - f_if_j(Z_{k}^{(n)})\middle\vert D }} \leq C\lambda_2^2e^{-C\lambda_2^2D} \epsilon_n. \]
	\end{lemma}
	
	Returning to the primary goal of this chapter, we finally can bound $\epsilon_{n+1}$ with the following proposition.
	\begin{proposition}\label{prop: contraction}
		Suppose that $\frac{2\sqrt{2} \max \pi_k^{3/2} \max \pi_k^{-3/2}}{\xi^3} \lambda_2 < \frac{\delta^2q^2}{8}$ and the conditions of \cref{cor: numChildCancel} hold. There exists another constant $C_1 = C_1(\xi, \delta, \pi ,q)$ such that if $\lambda_2^2d > C_1$, then for all $n$ sufficiently large, 
		\[ \epsilon_{n+1} \leq \frac{1}{2} \cdot \epsilon_n. \]
	\end{proposition}
	\begin{proof}
		Recall that 
		\begin{align*}
			\abs{(\mathcal{E}_{n+1})_{ij}} &= \abs{\E{ \sum_{k=1}^{D} \E{ f_if_j(Z_{k-1}^{(n)}) - f_if_j(Z_{k}^{(n)}) \middle\vert D}}} \leq \E{ \sum_{k=1}^{D} \abs{\E{ f_if_j(Z_{k-1}^{(n)}) - f_if_j(Z_{k}^{(n)}) \middle\vert D}}} \\
			&\leq \E{ \sum_{k=1}^{D} C\lambda_2^2e^{-C\lambda_2^2D}\epsilon_n} = \E{C\lambda_2^2D e^{-C\lambda_2^2D} \epsilon_n} = C\lambda_2^2 \epsilon_n\E{De^{-C\lambda_2^2D}}.
		\end{align*}
		Since $D \sim \text{Pois}(d)$, we can compute the expectation using the moment generating function. In particular, it is the first derivative of the moment generating function evaluated at $-C\lambda_2^2$. Since $\E{e^{tD}} = e^{d(e^t - 1)}$, we get that $\E{De^{tD}} = de^te^{d(e^t-1)}$. Since $e^t \leq 1 + C't$ for $t \in [-C, 0]$, we can evaluate at $t = -C\lambda_2^2$ and upper bound by $\E{De^{-C\lambda_2^2D}} \leq d(1 - C'C\lambda_2^2)e^{-C'C\lambda_2^2d} \leq de^{-C\lambda_2^2d}$. Substituting, we obtain
		\[ \abs{(\mathcal{E}_{n+1})_{ij}} \leq C\lambda_2^2de^{-C\lambda_2^2d} \epsilon_n. \]
		Now, for $\lambda_2^2d$ large enough in terms of $\pi, \xi, q, \delta$, we can ensure that the coefficient $C\lambda_2^2de^{-C\lambda_2^2d} < \frac{1}{2}$, which gives the desired result.
	\end{proof}
	The following corollary immediately follows.
	\begin{corollary}
		Suppose that $\frac{2\sqrt{2} \max \pi_k^{3/2} \max \pi_k^{-3/2}}{\xi^3} \lambda_2 < \frac{\delta^2q^2}{8}$, the conditions of \cref{cor: numChildCancel} hold, and $\lambda_2^2d > C_1$ as in \cref{prop: contraction}. Then 
		\[ \lim_{n \rightarrow \infty} \epsilon_n = \lim_{n\rightarrow \infty} \abs{\E{X_\rho^{(n)}(i) - \tilde X_\rho^{(n)}(i) \middle\vert \sigma_\rho = j}} = 0. \]
	\end{corollary}
	
	\section{The Algorithm}\label{section: algorithm}
	
	In the section, we provide the machinery needed to relate our previous results on the tree back to the stochastic block model. The machinery includes a structural coupling between local neighborhoods of the stochastic block model and our Galton-Watson tree, as well as a method to estimate the parameters on the tree given an instance of the block model. This will then allow us to outline our algorithm that recovers a provably asymptotically optimal fraction of the vertex label, which takes an initial black box reconstruction and amplifies it to an optimal reconstruction. 
	
	\subsection{Coupling Trees with the Stochastic Block Model}
	We provide the two coupling lemmas that will allow us to relate all of the previous work on trees to the stochastic block model itself. Much of the work was done in \cite{MNS:15} and remains largely the same, so we only provide a sketch of the argument and point at the steps at which an adjustment is needed. For a more rigorous presentation, refer to the original work by Mossel, Neeman, and Sly. 
	
	\begin{lemma}
		Let $R = \lfloor \frac{1}{10 \log(2(\max Q_{ij}))} \log n \rfloor$. For any fixed $v \in G$, there is a coupling between $(G, \sigma')$ and $(T, \sigma)$ such that $(B(v, R), \sigma_{B(v, R)}') = (T_R, \sigma_R)$ a.a.s.
	\end{lemma}
	
	\begin{proof}
		Here, $B(v, R)$ refers to the neighborhood of radius $R$ around $v$ in $G$, and $T_R$ refers to the first $R$ levels of the tree $T$. This notation is borrowed, and as mentioned above, we refer the reader to the proof of Proposition 4.2 in \cite{MNS:15}, and point out the differences in our setting. 
		
		For a vertex $v \in T$, we define $Y_v$ to be the number of children of $v$, and $Y_v^i$ to be the number of children whose label is $i$. Note that with these definitions, because of Poisson thinning
		\[ Y_v^i \sim \text{Pois}(dP_{\sigma_v i}) \]
		and that the pair $(T, \sigma)$ can be entirely reconstructed from $\sigma_\rho$  and the sequences $\{ (Y_u^i)_{u \in T} \}$. 
		
		For $G_R = B(v, R)$, we make similar definitions. Using the notation from \cite{MNS:15} that $V = V(G)$ and $V_R = V(G) \setminus V(B(v, R))$, we define $\{W^i\}$ to be the partition of $W \subset V$ into vertices with the corresponding label $i$. For any $v \in \partial G_R$, define $X_v$ to be the neighbors of $v$ in $V_R$. In this case, we have that 
		\[ X_v^i \sim \text{Binom}(\abs{V_R^i}, Q_{\sigma_v i}). \]
		The coupling relies on the fact that the Poisson and Binomial distributions with approximately the same expectations have asymptotically small total variation distance. However, the sequences of variables $X_v^i$ are not sufficient to reconstruct $G_R$. The reasons are because it is possible for two vertices $u, v$ to share children as we are no longer guaranteed to be working on a tree. The Lemmas 4.3 through 4.6 analyze this bad event, and all proofs go through the same as they do not depend on the number of communities present in the graph. Moreover, they do not rely on symmetry of communities either, but only that the edge probabilities are all $\Theta(\frac{1}{n})$. 
		
		Armed with these lemmas, we may now consider the conclusion to the proof. Here, the only argument which relies on the two communities is the event $\tilde \Omega$ which is defined as $\abs{\abs{V^+} - \abs{V^-}} \leq n^{3/4}$ and satisfies $\Prob{\tilde\Omega} \rightarrow 1$ exponentially fast. The essence of this condition is that the sizes of each community are close to their expected size. To achieve the same result, we define out event $\tilde \Omega$ to be $\abs{\abs{V^i} - n\pi_i} \leq n^{3/4}$ for all $i$. It remains to show that we still have $\Prob{\tilde \Omega} \rightarrow 1$ exponentially fast. Recall that $\abs{V^i} \sim \text{Binom}(n, \pi_i)$. By Hoeffding's Inequality, we know that
		\[ \Prob{\abs{\abs{V^i} - n\pi_i} \geq n^{3/4}} \leq 2\exp{-2n^{1/2}}. \]
		Union bounding over the $q$ communities, we see that indeed $\Prob{\tilde \Omega} \rightarrow 1$ exponentially fast as we wanted. Finally, in order to conclude instead of union bounding over 2 communities we union bound over $q$ communities. In either case, the union bound is over a quantity constant relative to $n$, so the proof concludes with the same result.
	\end{proof}
	\begin{lemma}
		For any fixed $v \in G$, there is a coupling between $(G, \tau')$ and $(T, \tau)$ such that $(B(v, R), \tau_{B(v, R)}') = (T_R, \tau_R)$ a.a.s. where $\tau'$ are the labels produced our black box.
	\end{lemma}
	
	\begin{proof}
		For this proof, we refer to the argument following Lemma 5.9 in \cite{MNS:16}. In the same spirit, we condition on $\sigma', B(v, R-1)$ and $G'$ and show that the conditional distribution of $\tau'$ is close to the distribution of $\tau$ conditioned on $\sigma$ and $T$. For any $u \in \partial B(v, R-1)$, we have that 
		\[ \abs{\{ wu \in E(G):  w \in G', \sigma'_w = i, \tau'_w = j \}} \sim \text{Binom}\left( \abs{V^i \cap W_v^j}, Q_{\sigma_u i}\right). \]
		We define $\Delta_{ij} = \frac{1}{n\pi_i} \abs{V^i \cap W_v^j}$. With this definition we have that $\abs{V^i \cap W_v^j} = \Delta_{ij} \cdot n\pi_i \pm O(n^{1/2})$, so by Lemma 4.6 in \cite{MNS:15}, the above distributions are at total variation distance at most $O(n^{-1/2})$ from $\text{Pois}(dP_{\sigma_u i}\Delta_{ij})$. Notice moreover that on the noisy tree $T$, we have for $u \in L_{R-1}$,
		\[ \abs{\{ wu \in E(T): \sigma_w = i, \tau_w = j \}} \sim \text{Pois}\left( dP_{\sigma_u i} \Delta_{ij} \right). \]
		In particular, the conditional distributions of $\tau$ and $\tau'$ on level $R$ are at total variation distance at most $O(n^{-1/2})$. Union bounding over the $O(n^{1/8})$ choices for $u$, we see that the two distributions are a.a.s the same, which gives the desired coupling. 
		
		We check that the $\Delta$ defined in this way does in fact satisfy the assumptions we specified in the main definitions. With sufficient accuracy, which we will have, we can ensure that the entries on the diagonal of $\Delta$ are close to 1, and the off-diagonal entries are close to 0. Specifically, by choosing $C$ to be large enough, we can guarantee that $\Delta_ii \geq 1 - \frac{1}{q^2}$ and $\Delta_{ij} \leq \frac{1}{q^2}$ for $i \neq j$. In this way, it is easy to see that out matrix $\Delta$ will be close enough to the identity that it will be invertible. More rigorously, we may appeal to the fact that the invertible matrices form an open subset of $\mathbb{R}^{q \times q}$. 
	\end{proof}
	
	\subsection{Parameter Estimation}
	We show that it is possible to accurately guess the parameters of our process on the tree from the provided graph. This will be essential to our algorithm, as in order to initialize the predictions on the boundary of the neighborhood, we calculate Bayesian priors that will depend on these parameters. To prove the result, we make use of the probabilistic method. 
	\begin{lemma}
		Given a set $U \subset V(G)$ such that $\abs{U} = \sqrt{n}$, we have at least $\Theta(n^{1/4})$ vertices with degree $k = \frac{1}{4} \cdot \frac{\log n}{\log \log n}$ a.a.s.
	\end{lemma}
	\begin{proof}
		First, we note that our degree $k$ satisfies $k^k \sim n^{1/4}$. Indeed, 
		\begin{align*}
			\log k^k &= k \log k = \frac{1}{4} \cdot \frac{\log n}{\log \log n} \cdot ( \log \log n - \log \log \log n - \log 4) \sim \frac{1}{4} \log n.
		\end{align*}
		Let $d$ be the average degree. We also have that
		\begin{align*}
			\Prob{\text{Bin}(n, \frac{d}{n}) = k} &= \binom{n}{k} \left( \frac{d}{n} \right)^k \left(1 - \frac{d}{n}\right)^{n-k} \geq \left( \frac{n}{k} \right)^k \left( \frac{d}{n} \right)^k e^{-d} = \frac{d^k}{k^k} e^{-d} \geq C n^{-1/4} e^{-d}.
		\end{align*}
		If we let $N$ be the number of vertices in $U$ that have degree $k$, then 
		\[ \E{N} \geq \sqrt{n} \cdot Cn^{-1/4}e^{-d} \geq \Theta(n^{1/4}).  \]
		Note that $N$ can be expressed as the sum of indicator variables $N = \sum_{u \in U} \mathbbm{1}(A_u)$ where $A_u$ is the event that $\deg(u) = k$. In particular, to analyze the distribution of $N$, we would like to understand 
		\[ \Delta = \sum_{u, v \in U} \text{Cov}(A_u, A_v). \]
		These two events are independent given the status of the edge between $u$ and $v$, so we expand the probability accordingly. 
		\begin{align*}
			\Delta &= \sum_{u, v \in U} \Prob{uv \in E(G)}\Prob{\text{Bin}(n-1, \frac{d}{n}) = k-1)}^2  + \Prob{uv \not \in E(G)}\Prob{\text{Bin}(n-1, \frac{d}{n}) = k)}^2 - \Prob{A_u}^2  \\
			&\approx \frac{C}{n} \cdot \Prob{\text{Bin}(n-1, \frac{d}{n}) = k-1}^2 - \frac{C}{n} \cdot \Prob{\text{Bin}(n, \frac{d}{n}) = k}^2 \approx \frac{C}{n}.
		\end{align*} 
		Since $N$ is the sum of indicators, we have the inequality that
		\[ \Var{N} \leq \E{N} + \sum_{u, v \in U} \text{Cov}(A_u, A_v). \]
		In particular,we have that
		\begin{align*}
			\Var{N} &\leq \Theta(n^{1/4})  + \binom{\sqrt{n}}{2} \cdot \frac{C}{n} = \Theta(n^{1/4}) + C = o(\E{N}^2).
		\end{align*}
		Thus, by an application of Chebyshev's inequality, we find that $N \sim \E{N}$ a.a.s.
	\end{proof}
	
	With this result, we can accurately estimate the noise matrix corresponding with our black box predictions. In particular, notice that with degrees tending to infinity, we can almost surely estimate the label of that corresponding vertex. Moreover, the number of such high degree vertices also tends to infinity, so we can accurately find at least one vertex of each community almost surely. Given a representative of high degree from each community, we can then simply estimate $\Prob{\sigma_v = j \middle\vert \sigma_u = i}$ by looking at the fraction of vertices adjacent to the representative from community $i$ whose label is $j$. This quantity can be written as an expression depending only on $\Delta_{ij}$ using Bayes' Rule, and so we can accurately solve for each entry $\Delta_{ij}$ using this process. 
 
\subsection{Black Box Algorithm (Sphere Comparison)}
As the black box algorithm to provide our initial noisy estimates, we will make minor adjustments to the following theorem from Abbe and Sandon \cite{AS:15}. 
\begin{theorem}\label{blackbox}
    For any $k \in \mathbb{Z}$, $p \in (0, 1)^k$ with $\abs{p} = 1$. and symmetric matrix $Q$ with no two rows equal, there exists $\epsilon(c)$ such that for all sufficiently large $c$, it is possible to detect communities in SBM$(n, p, cQ/n)$ with accuracy $1-e^{-\Omega(c)}$. 
\end{theorem}
We now describe the alterations we make to this algorithm for it to more suitably fit our needs. We would like for the partition to be proportioned as close to $\pi$ as possible. Notice that if there are at most $\epsilon n$ vertices classified incorrectly, to balance out the partition we re-assign at most $\epsilon n$ in total. In the worst case, all of these were originally labeled correctly, so we can guarantee that the partition balanced according to $\pi$ that has accuracy at least $1 - 2\epsilon$. 
	
	\begin{remark}
		The provided black box algorithm is given with conditions that are stronger than the ones we imposed to complete the proof. As such, given a new black box algorithm that achieves high accuracy under looser conditions, our results would imply an optimal algorithm under those conditions for free. Of course, this is only possible up to the conditions that we imposed, but there is room for easy improvement in this aspect. 
	\end{remark}
	
	\subsection{Optimal Algorithm}
	Finally, we can present the algorithm that produces the optimal partition of vertices. Our algorithm closely follows the same basic structure as the one presented in \cite{MNS:16}, and so the proof of its optimality is very similar as well, given our new results.
	
	\begin{algorithm}[H]\label{algorithm}
		\SetAlgoLined
		$R \leftarrow \lfloor \frac{1}{10 \log (2(\max Q_{ij}))} \log n \rfloor$ \\
		Take $U \subset V$ to be a random subset of size $\lfloor \sqrt{n} \rfloor$ \\
		$\{ W_*^i \} \leftarrow \emptyset$ \\
		$\{ W_{align}^i \} \leftarrow \text{BBPartition}(G)$ \\
		\For{$v \in V \setminus U$}{
			$\{ W_v^i \} \leftarrow \text{BBPartition}(G \setminus B(v, R-1))$ \\
			Monte Carlo estimate $\Delta$ using the high degree vertices in $U$ \\
			Permute $\{ W_v^i \} $ to align with $\{ W_{align}^i \}$ \\
			Define $\tau' \in [q]^{\partial B(v, R)}$ by $\tau'_u = i$ if $u \in W_v^i$ \\
			Add $v$ to $W_*^{\arg \max (\tilde X_{v, R}(\tau'))}$ \\
		}
		\For{ $v \in U$ } {
			Assign $v$ to one of the $W_*^i$ uniformly at random
		}
		\caption{Optimal Reconstruction Algorithm}
	\end{algorithm}
	
	With the proven theorems, to check correctness, we only need to check that we can align the various calls to the black box algorithm. Since we are guaranteed at least $1 - \epsilon$ correct vertices in each run, a correct alignment between runs will have at least $1 - 2\epsilon$ vertices in common. On the other hand, an incorrect pairing will have at most $2 \epsilon$ vertices in common. Going through each pair one by one, an alignment between calls takes at most $\binom{q}{2}$ such comparisons. With this alignment, we can run Algorithm \ref{algorithm}, and the couplings in the previous section show that this algorithm is in fact optimal in the case of $q$ communities in general, which is what we set out to achieve in this paper. 

\bibliographystyle{alpha}
\bibliography{SBM}

\end{document}